\documentclass{article}
\PassOptionsToPackage{numbers, compress}{natbib}
\usepackage[final]{nips_2017}

\usepackage[utf8]{inputenc} 
\usepackage[T1]{fontenc}    
\usepackage{hyperref}       
\usepackage{url}            
\usepackage{booktabs}       
\usepackage{amsfonts}       
\usepackage{nicefrac}       
\usepackage{microtype}      
\usepackage{graphicx} 
\usepackage{subfigure} 
\usepackage{bm}
\usepackage{natbib}
\usepackage{amsmath}
\usepackage{amsthm}
\usepackage{amssymb}

\theoremstyle{plain}
\newtheorem{thm}{\protect\theoremname}
\theoremstyle{definition}

\theoremstyle{plain}
\newtheorem{fact}[thm]{\protect\factname}
\theoremstyle{plain}
\newtheorem{lem}[thm]{\protect\lemmaname}

\newtheorem{remark}[thm]{\protect\remarkname}
\newtheorem{example}[thm]{\protect\examplename}

\providecommand{\factname}{Fact}
\providecommand{\lemmaname}{Lemma}
\providecommand{\problemname}{Problem}
\providecommand{\theoremname}{Theorem}
\providecommand{\remarkname}{Remark}
\providecommand{\examplename}{Example}

\newcommand{\Prob}{\mathbb{P}}
\newcommand{\bX}{\mathbf{X}}
\newcommand{\bY}{\mathbf{Y}}
\global\long\def\dotleq{\overset{\cdot}{\leq}}
\global\long\def\doteqdel#1{\overset{\cdot}{\underset{#1}{=}}}

\global\long\def\covind{\overset{d}{\to}}
\global\long\def\covinp{\overset{p}{\to}}
\global\long\def\bin{Bin}

\global\long\def\ind{\mathbb{\mathbb{I}}}
\global\long\def\E{\mathbb{E}}

\title{\Large{Minimax Optimality of Sign Test for Paired  Heterogeneous Data}}
\author{Martin J. Zhang
\\ Stanford University
\\ jinye@stanford.edu
\And Meisam Razaviyayn
\\ University of Southern California
\\ razaviya@usc.edu
  \And David Tse 
  \\ Stanford University
  \\ dntse@stanford.edu}
  
\begin{document}

%

%

%
%
%

\maketitle

\begin{abstract}
Comparing two groups under different conditions is ubiquitous in the biomedical sciences. In many cases, samples from the two groups can be naturally paired; for example a pair of samples may come from the same individual under the two conditions. However samples across different individuals may be highly heterogeneous. Traditional methods often ignore such heterogeneity by assuming the samples are identically distributed. In this work, we study the problem of comparing paired heterogeneous data by modeling the data as Gaussian distributed with different parameters across the samples. We show that in the minimax setting where we want to maximize the worst-case power, the sign test, which only uses the signs of the differences between the paired sample,  is optimal in the one-sided case and near optimal in the two-sided case. The superiority of the sign test over other popular tests for paired heterogeneous data is demonstrated using both synthetic data and a real-world RNA-Seq dataset.
\end{abstract}

\section{INTRODUCTION} \label{Intro}
A common form of scientific experimentation is the comparison of two groups. 
Suppose we collected $2n$ samples $\{X_i^A\}_{i=1}^n$, $\{X_i^B\}_{i=1}^n$ under two conditions $A$ and $B$.
The conditions may be sick v.s. healthy, pre- v.s. post- treatment, etc.
In the traditional homogeneous setting, samples within each group are assumed to be independently and identically distributed (i.i.d.), i.e. 
\begin{align*}
X_i^A \overset{\text{i.i.d.}}{\sim} \Prob_A,~~X_i^B \overset{\text{i.i.d.}}{\sim} \Prob_B,~~\forall~i=1,\cdots,n,
\end{align*}
where the distributions $\Prob_A$ and $\Prob_B$ typically come from some common distribution families like Gaussian or Poisson.
The goal is to infer whether there is a difference between the mean of the two groups, i.e. if $\E[ X_i^A ]\neq\E[ X_i^B]$.
One of the most commonly used test in this case is the two-sample t-test, which assumes the distributions $\Prob_A$ and $\Prob_B$ are Gaussian and is based on the t-statistic \cite{cressie1986use}.
However, in many real-world applications, the data are \textit{paired} and \textit{heterogeneous}. 
The paired structure means that for each $i$, the samples $X_i^A$, $X_i^B$ are similar due to some shared properties. 
The heterogeneity means that the data within the same group, $\{X_i^A\}$ or $\{X_i^B\}$, may be non-identically distributed. 
As a result, the paired differences $\{X_i^B-X_i^A\}$ may also be non-identically distributed.

Such paired heterogeneous data may occur in many scenarios. 
For example, in pre- v.s. post-treatment studies \cite{qian2013identification}, samples were taken before and after the treatment from the same individuals. 
Samples from the same person are similar and thus can be paired, while samples from different individuals may be very different due to individual-level heterogeneity. 
In another study, samples were obtained from the same person over a long period time to study viral infection disease (VID) \cite{chen2012personal}. 
Samples under different conditions (sick/healthy) can be paired if they are close to each other in time. 
As the person may change a lot over time, within-pair samples are more similar than within-group samples that are far from one another in time\footnote{One may argue that a time-series analysis is more appropriate \cite{aijo2014methods}. 
However, when we are not interested in the time-series pattern, the differential expression analysis by two-group comparison is still a valid method and has been used in various studies. 
Second, many RNA-Seq datasets including VID are noisy and have very few samples.
In those cases, the two-group comparison can produce more stable results compared to the time-series analysis.
}. See the following example.

\begin{figure}
\centering
\includegraphics[width=0.45\textwidth]{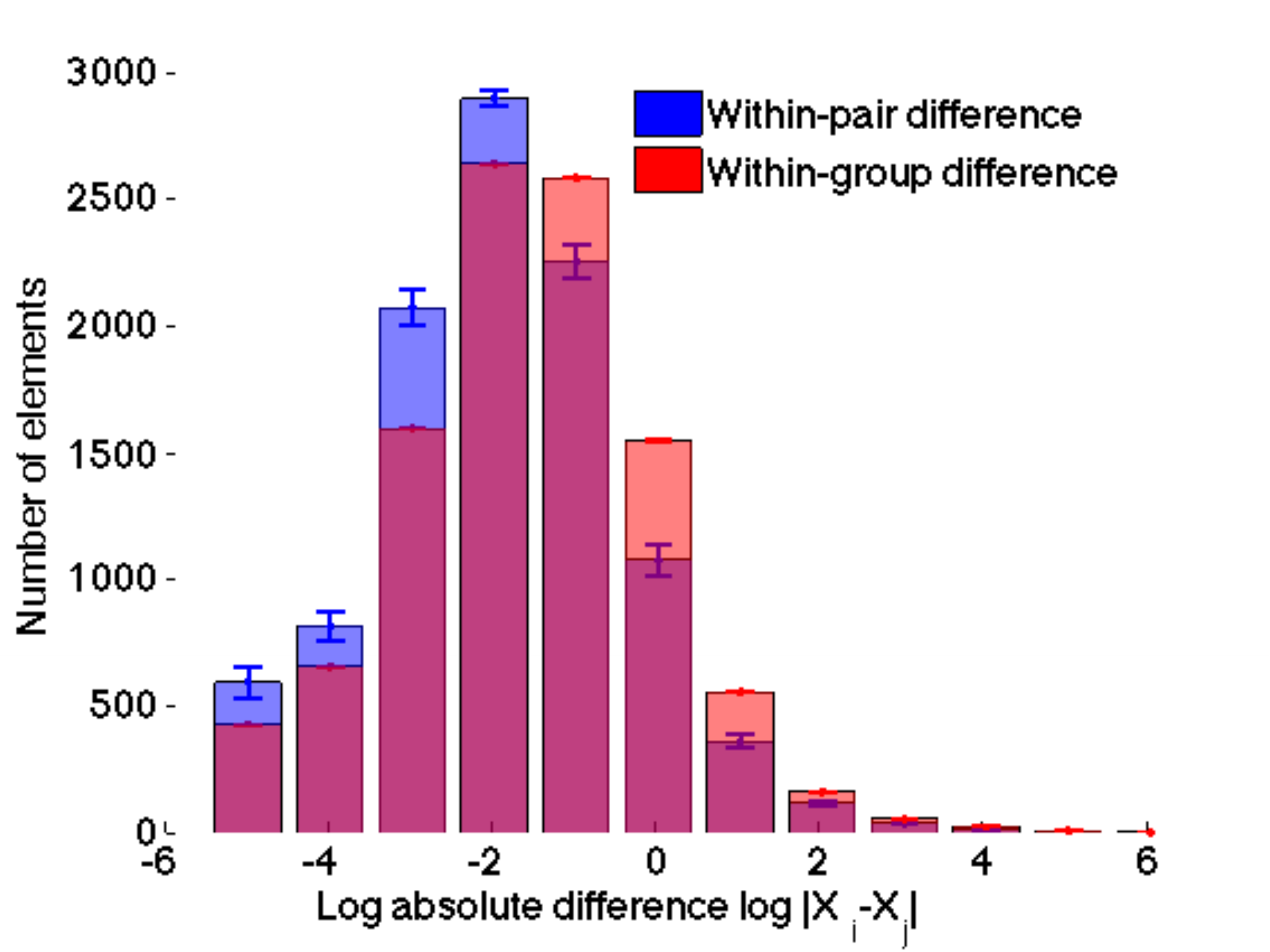}
\caption{Visualization of within-pair difference and within-group different for VID.\label{fig:EviDif}}
\end{figure}


\begin{example} (Heterogeneity in VID data)
We present a visualization on VID to illustrate the presence of the heterogeneity in data.
In this dataset, a sample $\mathbf{X}_i\in\mathbb{R}^{23,231}$ is a measurement of the gene expression level of $23,231$ genes, and samples are taken from one person over $1124$ days under two conditions, healthy and sick (see Fig. \ref{fig:VirInf}).
We match samples close to each other in time, one from each group (healthy/infected), as pairs. 
We plot the histogram of the pairwise difference of within-pair samples and of within-group samples (see the details in Supp. Sec. \ref{sec:vis}).
As can be seen in Figure~\ref{fig:EviDif}, the within-pair difference is systematically smaller than the within-group difference, indicating that the within-pair samples are more similar than within-group samples. 
\end{example}


A popular treatment to the paired heterogeneous data is to allow each data point to have a different distribution but assume that the paired differences $X_i^B - X_i^A$ are i.i.d.
If the paired differences are Gaussian, then the paired t-test is most powerful \cite{van2000asymptotic}.
The Gaussian assumption is quite reasonable and used in many applications, e.g. the RNA-Seq data \cite{law2014voom}.
However, given the heterogeneity across data pairs, it is hard to justify that the paired differences are actually identically distributed.
For example, in the pre- and post- treatment studies, different individuals may have different responses to the treatment, and hence the paired differences may not be i.i.d. 

In this work we keep the Gaussian assumption but allow the paired differences to be non-identically distributed.
To characterize the \textit{systematic} difference between the two group means, we assume the probability of increase/decrease from one group to the other is the same across all pairs. 
Specifically, we assume:
\begin{itemize}
\item Independently Gaussian (possibly non-identical):
\begin{align*}
X_i^A \sim \mathcal{N}(\nu_i^A, (\sigma^A_i)^2),~X_i^B \sim \mathcal{N}(\nu_i^B, (\sigma^B_i)^2),~\forall i.
\end{align*}
\item The ``tendency of shift'' is the same across all data pairs: $\Prob(X_i^B \geq  X_i^A) = \theta, \forall i$.
\end{itemize}
The tendency assumption is \textit{weaker} than the previous i.i.d. assumption and allows the paired differences to be non-identically distributed. 
In the Gaussian case, it implies a natural scaling for the paired differences $X_i^B-X_i^A$; their means being proportional to their standard deviations.
Even if the tendency assumption violated, the proposed test in this manuscript will still maintain a good power and be minimax optimal under certain conditions. 
See Remark \ref{rmk:vio_cond} for details.

Then the problem of interest is to test if $\theta=0.5$.
We seek robust tests that consistently produce high power under different levels of heterogeneity. A natural approach is to consider the minimax setting, where we fix the level of shift $\theta$ and maximize the worst-case power over all values of nuisance parameters. 

\textbf{Contributions.} 
The main contribution of the paper is to identify the optimal test for our minimax setting, which turns out to be the sign test \cite{lehmann2006testing}.
The sign test uses the number of times of the paired differences being positive as the summary statistics, i.e. $W = \sum_i \ind_{\{X_i^B-X_i^A>0\}}$. 
Our result shows that the sign test is maximin in the one-sided case where we want to test $\theta=0.5$ against $\theta > 0.5$. 
In the two-sided case where the alternative hypothesis becomes $\theta\neq 0.5$, we show that the worst-case power of any test can be upper bounded by that of the sign test plus a negligible additive term, implying the sign test is near optimal. 
In addition, we verify our theoretical analysis using both synthetic data and a real-world RNA-Seq dataset.

Let us explain our contributions within the context of prior art. 
Let the paired difference be $Y_i=X_i^B-X_i^A$. 
Prior to this work, it is known that if we restrict our attention to  $Y_i$'s only, which is natural given the paired structure, then the sign test is maximin over the class of all distributions where the differences $Y_i$'s are independent and the ``tendency of shift'' $\Prob(Y_i)$ is the same across all pairs \cite{lehmann2006testing}. 
In this case, the worst-case distribution is any distribution pair $\Prob_{H_0}$, $\Prob_{H_1}$ satisfying $\forall~i,y_+\geq 0,y_- < 0$,
\begin{align}\label{eq:worst_dist}
\frac{\Prob_{H_1}(Y_i=y_+ \vert Y_i \geq 0)}{\Prob_{H_0}(Y_i=y_+ \vert Y_i \geq 0)}= \frac{\Prob_{H_1}(Y_i=y_- \vert Y_i < 0)}{\Prob_{H_0}(Y_i=y_- \vert Y_i < 0)}=1.
\end{align}
The result is a direct consequence of plugging the above worst-case distribution in Theorem 8.1.1 in \citep{lehmann2006testing}.
The above argument considers a class of distributions so general that it may yield an overly pessimistic result. 
Indeed, the worst-case distribution $\Prob_{H_1}$ is very artificial in that it is continuous everywhere else but at $y=0$.
Hence, it is natural to restrict ourselves to a smaller and more natural class. 


According to empirical studies, the RNA-Seq data can be modeled as Gaussian random variables after variance-stabilization transformation \cite{law2014voom}.
In this work, therefore, we restrict ourselves to the family of  normal distributions.
As a result, the paired differences $Y_i$'s now have normal distributions.
Ideally, this distribution information should be properly utilized.
The question is that whether it leads to a test more powerful than distribution-free tests. 
Our result states that even this extra information does not help us to go beyond the sign test, indicating the importance of the sign information for robust testing. 
We also note that {\it our result is not a straight forward extension of existing results. In fact, it is not even clear that whether the sign test remains optimal after restricting to the Gaussian class because the Gaussian assumption excludes the worst-case distribution \eqref{eq:worst_dist} in the general class}.

In terms of the novelty of the proof, the proof techniques here are completely different from the older proof. 
Theorem 8.1.1 \cite{lehmann2006testing} cannot be applied here because it is extremely difficult to find the worst-case distribution in the Gaussian case (See the discussion after Theorem \ref{thm: opt_signtest}). 
{\it In fact, our conjecture is that the worst-case distribution does not even exist.}
We used a different strategy in our proof: we first show that the sign test is maximin among the family of “simple tests”. Then we show that “simple tests” can approximate the Borel measurable tests arbitrarily well. This is inspired by the widely used techniques in measure theory. However, we made two changes here: first, the notion of “simple tests”  is tailored to fit the location-scale invariance property, different from simple functions in measure theory; second, the approximation is in terms of the testing performance (size and power), rather than some function norms.
From a theoretical point of view, our result fills in a missing piece in the minimax analysis (of the adaptivity of sign test to the Gaussian family) in the classical statistical literature. 

{\bf Related works.}
This work has a very classical flavor. 
Some related topics are testing within-group heterogeneity \cite{davison1992treatment}, robust tests for paired data \cite{grambsch1994simple}, and rank-based tests \cite{lehmann2006nonparametrics}.
In the literature, paired t-test and the Wilcoxon signed-rank test are compared most often to the sign test \cite{van2000asymptotic}.
The paired t-test assumes the paired differences $X_i^B - X_i^A$ are i.i.d. Gaussian and uses the t-statistic \eqref{eq:stats_pairt}. 
Let $S_i$ be the sign of $X_i^B - X_i^A$ and $R_i$ be the rank of $\vert X_i^B - X_i^A \vert$ among all pairs. The Wilcoxon test uses the sign-rank statistic \eqref{eq:stats_wil}.
\begin{align} \label{eq:stats_pairt}
&\text{Paired t-test}: T = \sqrt{n} \frac{\text{mean}(X_i^B - X_i^A)}{\text{std}(X_i^B - X_i^A)}, \\
 &\text{Wilcoxon test}: U=\sum _{{i=1}}^{n}S_i R_i. \label{eq:stats_wil}
\end{align}
Most of the existing results for the above  tests are based on the i.i.d. (homogeneous) scenario. For example, when the differences $X_i^B - X_i^A$ are i.i.d. Gaussian, the paired t-test is known to be most powerful and the relative efficiency of the sign test v.s. the paired t-test is $2/\pi$ (Table 14.1,\cite{van2000asymptotic}), and $3/\pi$ for Wilcoxon test (an extension of the former). Results are rare on the heterogeneous case, especially under the minimax setting.

The motivating application for the present work is RNA-Seq experiments.
In RNA-Seq experiments, the gene expression level of people under different conditions are measured and the task is to identify genes differentially expressed under the two conditions. 
In related works, within-group heterogeneity is usually modeled by assuming some prior distribution on the expression level, e.g. gamma distribution \cite{chung2013differential, hardcastle2013empirical}. 
The paired structure is modeled by the design matrix in the generalized linear model \cite{law2014voom, love2014moderated}, or by assigning same expression level parameters to samples in the same pair \cite{chung2013differential, hardcastle2013empirical}. 
All above methods assume some complex models for the data, e.g. Bayesian hierarchical model in \cite{chung2013differential} or some mean-variance function shared across genes in \cite{robinson2010edger, love2014moderated}.
This leads to the lack of thorough theoretical understanding and, consequently, difficulty in establishing theoretical guarantees. 

\textit{``Those methods treat the estimated parameters as if they were known parameters, without allowing for the uncertainty of estimation, and this leads to statistical tests that are overly liberal in some situations'' \cite{soneson2013comparison}}. 

In fact, no theoretical result is yet available for the RNA-Seq data analysis \cite{law2014voom}. 
On the contrary, the sign test is theoretically justified in this manuscript. As shown in the experiments, it can be easily applied to the RNA-Seq data after simple normalization. Moreover, despite its simplicity, it yields reasonable results compared to other much more complex methods.

The rest of the paper is organized as follows. 
After the problem formulation in Sec. \ref{sec:PblmFor}, we prove the optimality of the sign test in Sec. \ref{sec:OptSign}, followed by a theoretical comparison of the sign test with the paired t-test in SubSec. \ref{sec:comp}.
Finally, we present numerical experiments on both synthetic data and real-world data in Sec. \ref{sec:Exp}.
We postpone the proofs to the supplementary materials.

\section{PROBLEM FORMULATION} \label{sec:PblmFor}
Consider $n$ paired data points $\{X_i^A,X_i^B \}_{i=1}^n$, where $(X_i^A,X_i^B)$ denotes the $i$-th sample pair in groups $A$ and $B$. 
Our goal is to detect whether there is a systematic difference between samples in two groups. 
We assume that 1. the samples are independently and normally distributed; 2. the ``tendency of shift'' is the same across all sample pairs.
Let $[n]$ denote the set $\{1,2,\cdots,n\}$. Then mathematically the above assumptions can be written as
\begin{align*}
& X_{i}^{A}\sim {\cal N}(\nu_{i}^{A},(\sigma_{i}^{A})^{2}), ~~
X_{i}^{B}\sim{\cal {\cal N}}(\nu_{i}^{B},(\sigma_{i}^{B})^{2}),~~\forall i \in [n], \\
& \Prob(X_{i}^{B}\geq X_{i}^{A})= \Prob(X_{j}^{B}\geq X_{j}^{A}) \triangleq \theta, \;\forall i,j \in [n],
\end{align*}
The range of parameters are
\begin{align*}
\{\nu_{i}^{A}, \nu_{i}^{B} \in \mathbb{R},~\sigma_{i}^{A}, \sigma_{i}^{B} \in \mathbb{R}_{\geq 0},~\textrm{s.t.}~\Prob(X_{i}^{B}\geq X_{i}^{A})=\theta\},
\end{align*} 
where $\mathbb{R}_{\geq 0}$ denotes the set of non-negative real numbers.
Clearly, $\theta = 0.5$ means that there is no systematic shift from group $A$ to group $B$, and $\theta \neq 0.5$ indicates that such systematic difference exists. 
Hence, our null hypothesis is $\theta = 0.5$.
For the alternative hypothesis, if we have some prior knowledge on the shifting direction, we can test a one-sided alternative $\theta>0.5$. 
Otherwise, a two-sided alternative, $\theta \neq 0.5$, is appropriate. 

Let us represent our statistical model in a more tractable way. Let $\Phi(\cdot)$ be the cumulative density function of the standard normal distribution. Since $\forall i$, $\Prob(X_{i}^{B} \geq X_{i}^{A})=\Phi(-\frac{\nu_{i}^{B}-\nu_{i}^{A}}{\sqrt{(\sigma_{i}^{B})^{2}+(\sigma_{i}^{A})^{2}}})=\theta$, we obtain $\forall i,j$, $
\frac{\nu_{i}^{B}-\nu_{i}^{A}}{\sqrt{(\sigma_{i}^{B})^{2}+(\sigma_{i}^{A})^{2}}} = \frac{\nu_{j}^{B}-\nu_{j}^{A}}{\sqrt{(\sigma_{j}^{B})^{2}+(\sigma_{j}^{A})^{2}}}$.
Let $\delta \triangleq \frac{\nu_{i}^{B}-\nu_{i}^{A}}{\sqrt{(\sigma_{i}^{B})^{2}+(\sigma_{i}^{A})^{2}}}$. Clearly, $\theta = \Phi(-\delta)$, and by defining  $\mu_{i} \triangleq \sqrt{(\sigma_{i}^{A})^{2}+(\sigma_{i}^{B})^{2}}$, $\rho_{i}\triangleq \frac{(\sigma_{i}^{A})^{2}}{(\sigma_{i}^{A})^{2}+(\sigma_{i}^{B})^{2}}$,
and $\nu_{i}\triangleq \nu_{i}^{A}$, the above model becomes 
\begin{equation}
\label{eq:SampleModel}
\begin{array}{ll}
&X_{i}^{A}\sim{\cal {\cal N}}(\nu_{i},\rho_{i}\mu_{i}^{2}), \\
& X_{i}^{B}\sim{\cal {\cal N}}(\nu_{i}+\delta\mu_{i},(1-\rho_{i})\mu_{i}^{2}),~ i \in [n],
\end{array}
\end{equation}
where $\nu_{i} \in \mathbb{R}$, $\mu_{i} \in \mathbb{R}_{>0}$, $\rho_{i} \in [0,1]$. Here  $\mathbb{R}_{>0}$ is the set of positive real numbers and the tendency assumption prevents $\mu_i$ to be $0$.
It is not hard to see that in this equivalent representation, we test the null hypothesis $\mathcal{H}_0: \delta = 0$ against the alternative hypothesis $\mathcal{H}_1: \delta > 0$ (one-sided) or $\mathcal{H}_1: \delta \neq 0$ (two-sided). 
To quantify the size and the power, we let $\delta$ have some fixed unknown magnitude. For the two-sided case, it may be either positive or negative with the sign $s_\delta$.
Since we have no knowledge about the nuisance parameters  $\{\nu_{i},\mu_{i},\rho_{i},s_\delta\}_{i=1}^n$, a natural formulation is to look for the maximin test that maximizes the worst-case power over all possible values of the nuisance parameters. 

Let the data vectors $ \bX^A\triangleq \{X_i^A\}_{i=1}^n$ and $\bX^B \triangleq \{X_i^B\}_{i=1}^n$ be data points obtained by model~\eqref{eq:SampleModel}. We use $\phi(\mathbf{X}^{A},\mathbf{X}^{B}):\mathbb{R}^{n}\times\mathbb{R}^{n}\mapsto[0,1]$ to denote a test that rejects the null hypothesis with probability $\phi(\bX^A, \bX^B)$ when the data are $\bX^A, \bX^B$. Given the nuisance parameters  $\gamma \triangleq\{\nu_{i},\mu_{i},\rho_{i},s_\delta\}_{i=1}^n$, the size is given by $\E_{\Prob_0(\gamma)} [\phi(\bX^A,\bX^B)]$, where the expectation is taken with respect to the null distribution $\Prob_0(\gamma)$ with the given nuisance parameters $\gamma$. Similarly, the power of the test is given by $\E_{\Prob_1(\gamma)} [\phi(\bX^A,\bX^B)]$.
We call a test $\phi^*(\cdot,\cdot)$, a level-$\alpha$ maximin test if for any other test $\phi$,
\begin{align*}
&\inf_{\gamma} \E_{\Prob_1(\gamma)} [\phi^*(\bX^A,\bX^B)] \geq \inf_{\gamma} \;\E_{\Prob_1(\gamma)} [\phi(\bX^A,\bX^B)],\\&
\sup_\gamma \E_{\Prob_0(\gamma)} [\phi^*(\bX^A,\bX^B)] \leq \alpha.
\end{align*}
In other words, $\phi^*$ has the best worst-case power among all tests with size smaller than $\alpha$ over all values of the nuisance parameters $\gamma$. 
Equivalently, the problem can be stated as
\begin{align} \label{eq:pblmminimax}
&\phi^* \in \arg \max_{\phi} \inf_\gamma \E_{\Prob_1(\gamma)} [\phi(\bX^A,\bX^B)], \nonumber\\ & {\rm s.t.}\quad  \sup_{\gamma} \E_{\Prob_0(\gamma)} [\phi(\bX^A,\bX^B)] \leq \alpha .
\end{align}
For the sake of notational simplicity, we abbreviate the expressions as follows.
Given two vectors $\mathbf{a} $ and $\mathbf{b}$, let $\mathbf{a} \circ \mathbf{b}$ be the vector of the same dimension that contains the element-wise product of $\mathbf{a}$ and $\mathbf{b}$. 
Also by the inequalities $\mathbf{a}>\mathbf{b}$ we refer to element-wise comparison, i.e $a_i>b_i, \forall i$.
For two tests $\phi$ and $\psi$, we use $\phi\doteqdel{\epsilon}\psi$ to denote that the two tests have $\epsilon$-similar performance, 
\begin{align}\label{eq:test_apprx}
& \vert\inf_{\gamma}\mathbb{E}_{\Prob_1(\gamma)}[\phi(\bX^A,\bX^B)]-\inf_{\gamma}\mathbb{E}_{\Prob_1(\gamma)}[\psi(\bX^A,\bX^B)]\vert\leq\epsilon \nonumber\\
&\vert\sup_{\gamma}\mathbb{E}_{\Prob_0(\gamma)}[\phi(\bX^A,\bX^B)]-\sup_{\gamma}\mathbb{E}_{\Prob_0(\gamma)}[\psi(\bX^A,\bX^B)]\vert\leq\epsilon.
\end{align}
As a natural extension, $\phi\doteq\psi$ means that the two tests have the same performance, i.e., $\phi\doteqdel 0\psi$. Similarly, $\phi\dotleq\psi$ means $\phi$ has no better performance as $\psi$; in other words,
\begin{align*}
&\inf_{\gamma}\mathbb{E}_{\Prob_1(\gamma)}[\phi(\bX^A,\bX^B)] \leq \inf_{\gamma}\mathbb{E}_{\Prob_1(\gamma)}[\psi(\bX^A,\bX^B)]\\
&\sup_{\gamma}\mathbb{E}_{\Prob_0(\gamma)}[\phi(\bX^A,\bX^B)]\geq\sup_{\gamma}\mathbb{E}_{\Prob_0(\gamma)}[\psi(\bX^A,\bX^B)].
\end{align*}
Finally, we will use $\E_0$ and $\E_1$ as shorthand representations for $\E_{\Prob_0(\gamma)}$ and $\E_{\Prob_1(\gamma)}$, respectively. 

\section{OPTIMALITY OF SIGN TEST\label{sec:OptSign}} 
Let us start by showing the location-scale invariance of the minimax problem \eqref{eq:pblmminimax}. 
Generally speaking, the worst-case performance of a test $\phi(\cdot,\cdot)$ does not change under any shifting or scaling of the input:
\begin{fact}
\label{fact:scl_inv}
Let $\phi(\mathbf{X}^{A},\mathbf{X}^{B}):\mathbb{R}^{n}\times\mathbb{R}^{n}\mapsto[0,1]$
be an arbitrary test. For any $\mathbf{a}\in\mathbb{R}^{n}$, $\mathbf{b}\in\mathbb{R}_{>0}^{n}$,
let $\psi(\mathbf{X}^{A},\mathbf{X}^{B})=\phi(\mathbf{a}+\mathbf{b}\circ\mathbf{X}^{A},\mathbf{a}+\mathbf{b}\circ\mathbf{X}^{B})$.
Then, $\phi\doteq\psi$. 
\end{fact}
See Supp. Sec \ref{sec:pf_fact} for the proof.
Fact \ref{fact:scl_inv} leads to two conjectures.
First, any data shifting does not affect the performance, suggesting that the absolute offset may be redundant and we should only focus on the relative difference $\mathbf{Y}\triangleq \mathbf{X}^{B}-\mathbf{X}^{A}$.
Second, any data scaling does not affect the performance either. This suggests that the magnitude of the data, $\vert\mathbf{Y}\vert$, may also be redundant. 
Then, intuitively what remains, namely the signs of the difference $S_{i} \triangleq \text{sgn}(Y_{i})$, is the actual informative part for the minimax problem. 

For the most powerful test using only the sign information, a sufficient statistic is the number of positive signs $W=\sum_{i=1}^{n}\ind_{\{Y_{i}>0\}}$. Clearly,
under the null distribution,  $W\sim\bin(n,0.5)$. 
Then the sign test $\phi^S(\bX^A,\bX^B)$ can be written as 
\begin{align}
&\text{one-sided:}~\ind_{\{W>c_1\}} + p_1\ind_{\{W=c_1\}}, \label{eq:SgnTst_oneside}\\ &
\text{two-sided:}~\ind_{\{\vert W- n/2 \vert >c_2\}} + p_2\ind_{\{\vert W- n/2 \vert = c_2\}}, \label{eq:SgnTst_twoside}
\end{align}
where $c_1$, $c_2$, $p_1$, $p_2$ are some constants calculated according to $W\sim\bin(n,0.5)$ and the level $\alpha$.
Note that $p_1,p_2\in[0,1]$.
In addition, the distribution of $W$ and the power of the sign test does not depend on the values of the nuisance parameters.
We next prove the optimality of the sign test in the one-sided case and near optimality in the two-sided case. 

\subsection{One-sided Case}
In the one-sided case, let us assume $\theta\geq 0.5$ ($\delta \geq 0$). 
Then the nuisance parameters are just $\gamma \triangleq\{\nu_{i},\mu_{i},\rho_{i}\}_{i=1}^n$. 
Our main result below confirms that the sign test is indeed maximin in the sense of  \eqref{eq:pblmminimax}.
\begin{thm}
\label{thm: opt_signtest} 
Let $\mathcal{B}_{n\times n}$ be the class of Borel measurable functions that maps $\mathbb{R}^{n}\times\mathbb{R}^{n}$ to $[0,1]$. Then the one-sided sign test, as given in \eqref{eq:SgnTst_oneside}, is maximin among all tests $\phi(\mathbf{X}^{A},\mathbf{X}^{B})\in\mathcal{B}_{n\times n}$.
\end{thm}

Recall from Section \ref{sec:PblmFor} that our statistical model is
\begin{align*}
X_i^A \sim \mathcal{N}(\nu_i, \rho_i \mu_i^2),~~~~X_i^B \sim \mathcal{N}(\nu_i+\delta\mu_i, (1-\rho_i) \mu_i^2).
\end{align*}
Paired tests literatures suggest us to look at only the difference $Y_i=X_i^B - X_i^A \sim \mathcal{N}(\delta \mu_i, \mu_i^2)$ to get rid of the nuisance parameters $\nu_i$'s and $\rho_i$'s, which is not surprising. 
Here let us take the case $n=1$ as an example to give some high-level intuitions why we can further reduce the sufficient statistics from $Y_i$'s to $S_i$'s.  

According to Theorem 8.1.1 and the Neyman-Pearson lemma in \citet{lehmann2006testing}, to show the sign test is maximin, it suffices to find a prior on $\mu_1$ where $S_1$ is a sufficient statistic.
However, a careful inspection of the proof of Lemma \ref{lem:stepopt}, especially on \eqref{eq:stepopt_pf_4}, reveals that there is no such single prior on $\mu_1$ for which $S_1$ is a sufficient statistic.
However, there is in fact a sequence of priors on $\mu_1$ such that fixing the observation $Y_1$, $S_1$ is asymptotically a sufficient statistic. 
For $k=1,2,\cdots$, consider the sequence of priors $g_k(\mu_1) = c_k/ \mu_1$ for $\mu_1 \in (1/k, k)$ and some normalizing constant $c_k$. 

Let $f_0(\cdot;\mu)$, $f_1(\cdot;\mu)$ be the densities of $\mathcal{N}(0,\mu^2)$, $\mathcal{N}(\delta \mu,\mu^2)$ respectively
and let $f_0$, $f_1$ be that of $\mathcal{N}(0,1)$, $\mathcal{N}(\delta,1)$.
A direct calculation (by change of variable $\mu'=\frac{Y_1}{\mu}$) shows that as $k\to \infty$, the likelihood ratio 
\begin{align*}
&\frac{f(Y_1|H_1)}{f(Y_1|H_0)} = \frac{\int_{\frac{1}{k}}^{k} f_1(Y_1;\mu) \frac{1}{\mu} d \mu}{\int_{\frac{1}{k}}^{k} f_0(Y_1;\mu) \frac{1}{\mu} d \mu}  = \frac{\int_{\frac{1}{k}}^{k} f_1(\frac{Y_1}{\mu}) \frac{1}{\mu^2} d \mu}{\int_{\frac{1}{k}}^{k} f_0(\frac{Y_1}{\mu}) \frac{1}{\mu^2} d \mu} \overset{\mu'=\frac{Y_1}{\mu}}{=} \\&  \frac{\int_{\frac{Y_1}{k}}^{kY_1} f_1(\mu') d \mu'}{\int_{\frac{Y_1}{k}}^{kY_1} f_0(\mu') d \mu'} \overset{k\to \infty}{\to}
2 \left[\theta \ind_{\{S_1=+\}} + (1-\theta) \ind_{\{S_1=-\}}\right].
\end{align*}
Then asymptotically, the likelihood ratio of the observation  $f(Y_1|H_1)/f(Y_1|H_0)$ depends only on $S_1$, implying that $S_1$ is indeed asymptotically a sufficient statistic.

We note the above argument serves \textit{only} as a high-level intuition and is by no means rigorous. 
More specifically, here we only show that for each $Y_1$, $S_1$ is asymptotically a sufficient statistic for the sequence of priors $g_1, g_2,\cdots$, which is essentially a point-wise convergence result. 
However, a uniform convergence is needed to actually prove the final result. 
The rigorous proof is as follows.

\begin{proof}
(Proof sketch of Theorem \ref{thm: opt_signtest}) 
The idea is to show any test in $\mathcal{B}_{n\times n}$ performs no better than the sign test. The proof is composed of three main steps. We state the lemmas being used right after the corresponding steps, and relegate their proofs to the supplementary material. 
We recall that  inequalities on vectors are element-wise, e.g. $\mathbf{a} \geq \mathbf{b}$ means $a_i\geq b_i$ for all $i$.

\textbf{Step 1:} 
Define the set of Borel measurable tests:
\begin{align*}
\mathcal{B}_{n}=\{f:\mathbb{R}^{n}\mapsto[0,1],\;f\;is\;Borel\;measurable\}.
\end{align*}
Lemma \ref{lem:two_to_one} implies that it suffices to show that the sign test is maximin among all tests $\phi(\bY)\in\mathcal{B}_n$.  
\begin{lem}
\label{lem:two_to_one}
For any test $\phi(\mathbf{X}^{A},\mathbf{X}^{B})\in{\cal B}_{n\times n}$,
there exists a Borel measurable test $\psi(\bY) \in \mathcal{B}_n$, such that  $\phi\dotleq\psi$. Moreover, if $\phi$ is symmetric, then $\psi$ is also symmetric\footnote{The symmetry is used for proving Theorem \ref{thrm:opt_two_side}.}. 
\end{lem} 

\textbf{Step 2:} 
We show the sign test is maximin over the set of ``simple test'' $\mathcal{S}\subset \mathcal{B}_n$, which is defined as follows. 
Let $\mathcal{O} = \{-,+\}^n$ be the set of $2^n$ orthants in $\mathbb{R}^n$ and let $\mathbf{o}=(o_1,\cdots,o_n) \in \mathcal{O}$. 
Consider any $\omega>0$.
For any $b\in\mathbb{Z}$, let the $1$-D intervals be $I_{b}^{+}=((1+\omega)^{b},(1+\omega)^{b+1}]$ and $I_{b}^{-}=[-(1+\omega)^{b+1},-(1+\omega)^{b})$.
Define the $n$-D box, specified by the orthant index $\mathbf{o}$ and the interval index $\mathbf{b}=(b_1,\cdots,b_n)$, as $I_{\mathbf{b}}^\mathbf{o}=I^{o_1}_{b_1}\times \cdots \times I^{o_n}_{b_n}$. 
Then define the set of simple tests, denoted by $\mathcal{S}(\omega)$, to be the test that are piece-wise constant on the boxes $I_{\mathbf{b}}^\mathbf{o}$'s: 
\begin{equation}\label{eq:defSomega}
\begin{split}
\mathcal{S}(\omega)=&\{\phi:\phi= \sum_{\mathbf{o}\in\mathcal{O}} \sum_{-\infty<\mathbf{b}<\infty}\phi_{\mathbf{b}}^{\mathbf{o}}\ind_{I_{\mathbf{b}}^{\mathbf{o}}} +\phi_{0}\ind_{\{0\}}, \\&for~some~0\leq\phi_{\mathbf{b}}^{\mathbf{o}}\leq 1,0\leq\phi_0\leq 1 \},
\end{split}
\end{equation}
and let ${\cal S}=\bigcup_{\omega>0}{\cal S}(\omega)$. 
By Lemma \ref{lem:stepopt}, the sign test is maximin among all tests in ${\cal S}$. 

\begin{lem} 
\label{lem:stepopt}
The one-sided sign test \eqref{eq:SgnTst_oneside} is maximin among all $\alpha$-level tests in ${\cal S}$. 
\end{lem}

\textbf{Step 3:} We show that $\mathcal{S}$ approximates $\mathcal{B}_n$ arbitrarily well in terms of testing performance as defined in \eqref{eq:test_apprx}, and hence establish the optimally of the one-sided sign test in $\mathcal{B}_n$.
Specifically, by Lemma \ref{lem:sim_apr}, $\forall~\phi\in{\cal B}_n,~\epsilon>0$, $\exists~\psi\in{\cal S}$, s.t. $\phi\doteqdel{\epsilon}\psi$. 
Letting $\epsilon\downarrow0$ we have that $\phi^S$ is maximin among all tests in ${\cal B}_n$, concluding the proof. 
\begin{lem}
\label{lem:sim_apr} Let  ${\cal B}_n$ be the set of Borel measurable functions $f:\mathbb{R}^n\mapsto[0,1]$.
For any $\phi(\mathbf{Y})\in{\cal B}_n$ and any $\epsilon>0$, there exists a measurable function 
$\psi\in{\cal S}$ such that $\phi\doteqdel{\epsilon}\psi$. 
\end{lem}
\end{proof}

\begin{remark} \label{rmk:vio_cond} 
Under the alternative distribution, the testing statistics $W$ will follow a binomial distribution $\text{Bin}(n,\theta)$.
If the tendency assumption is violated, i.e. each pair has a different $\theta_i$, it will instead follow a Poisson binomial distribution with parameter $(\theta_1,\cdots,\theta_n)$, which has a tail property similar to that of the binomial distribution. 
Hence the sign test will still maintain a good power.

Moreover, when $\theta_i$'s are different, one can consider a minimax setting over $\theta_i$'s by defining the nuisance parameters to be $\{\mu_i,\rho_i,\nu_i\} \cup \{\theta_i: \theta_i \geq \theta_0\}$ for some $\theta_0>0.5$. 
It is not hard to see that the one-sided sign test is still maximin in this case.
Specifically, the power will increase with the increase of any $\theta_i$, and hence the worst-case is when $\theta_i=\theta_0~ \forall i$, which reduces to the setting where $\theta_i$'s are same. 
\end{remark}

\subsection{ Two-sided Case}
Now we extend our result to the two-sided case, where we want to test $\theta=0$ v.s. $\theta\neq 0$. 
Recall that in this case, we can no longer assume $\delta\geq 0$ for the distribution in (\ref{eq:SampleModel}).
So we modify the formulation in Section \ref{sec:PblmFor} by letting $s_\delta\in\{-1,1\}$ to be the sign of $\delta$, and letting the nuisance parameter be $\gamma = \{\nu_{i},\mu_{i},\rho_{i},s_\delta\}_{i=1}^n$. 
We fix the magnitude $\vert \delta \vert$ and consider the maximin problem \eqref{eq:pblmminimax}.
Without loss of generality assume $\alpha < 0.5$. 
Let $\tilde{\phi}^S$ be the $\frac{\alpha}{2}$-level one-sided sign test. 
The $\alpha$-level two-sided sign test $\phi^S$ can be written as 
\begin{equation} \label{def:two_side_sign}
\phi^S=\tilde{\phi}^S(\bY)+\tilde{\phi}^S(-\bY)
\end{equation}
The following theorem shows in the two-sided case, the sign test is near optimal. See Supp. Subsec. \ref{pf:opt_two_side} for the proof.

\begin{thm} (Two-sided case) \label{thrm:opt_two_side}
Let $\mathcal{B}_{n\times n}$ be the class of Borel measurable functions that maps $\mathbb{R}^{n}\times\mathbb{R}^{n}$ to $[0,1]$, and let $\phi^S$ be the two-sided sign test as defined in \eqref{def:two_side_sign}. For any $\alpha$-level test $\phi \in \mathcal{B}_{n\times n}$, the worst-case power satisfies
\begin{equation*}
\inf_\gamma \E_1[\phi] \leq \inf_\gamma \E_1[\phi^S] +  \frac{\alpha}{2}  \exp(-\frac{n \delta^2}{2}). 
\end{equation*}
If $\alpha=0.05$ and $\delta=\frac{3}{\sqrt{n}}$, the additive term is $2.7e$-$4$, almost negligible. 
\end{thm}

\begin{remark}
The proofs in both cases mainly use two properties of the Gaussian distribution. The first is the location-scale invariance, i.e., if we scale or shift the data points, the distribution still lies within the family of interest.  This is used in the proof of Fact \ref{fact:scl_inv}, Lemma \ref{lem:two_to_one}, and Lemma \ref{lem:stepopt}. The second is the sub-Gaussian tail property.
This is used in the proof of Lemma \ref{lem:sim_apr} and Theorem \ref{thrm:opt_two_side}. 
Specifically, if we have a heavy-tailed distribution but the tail probability still vanishes, Lemma \ref{lem:sim_apr} will still hold, while the exponential term on the RHS of Theorem \ref{thrm:opt_two_side} will become a term with a slower vanishing speed that depends on the tail property of the distribution under consideration.

To generalize to result, in order for the sign test to be maximin in some family of distribution, the family needs to have location-scale invariance and sub-Gaussian tail property. Our conjecture is that these two are also sufficient for the minimaxity of the sign test.
\end{remark}

\begin{figure*}
\centering
\subfigure[Power versus Magnitude. ]{\includegraphics[width=0.3\textwidth]{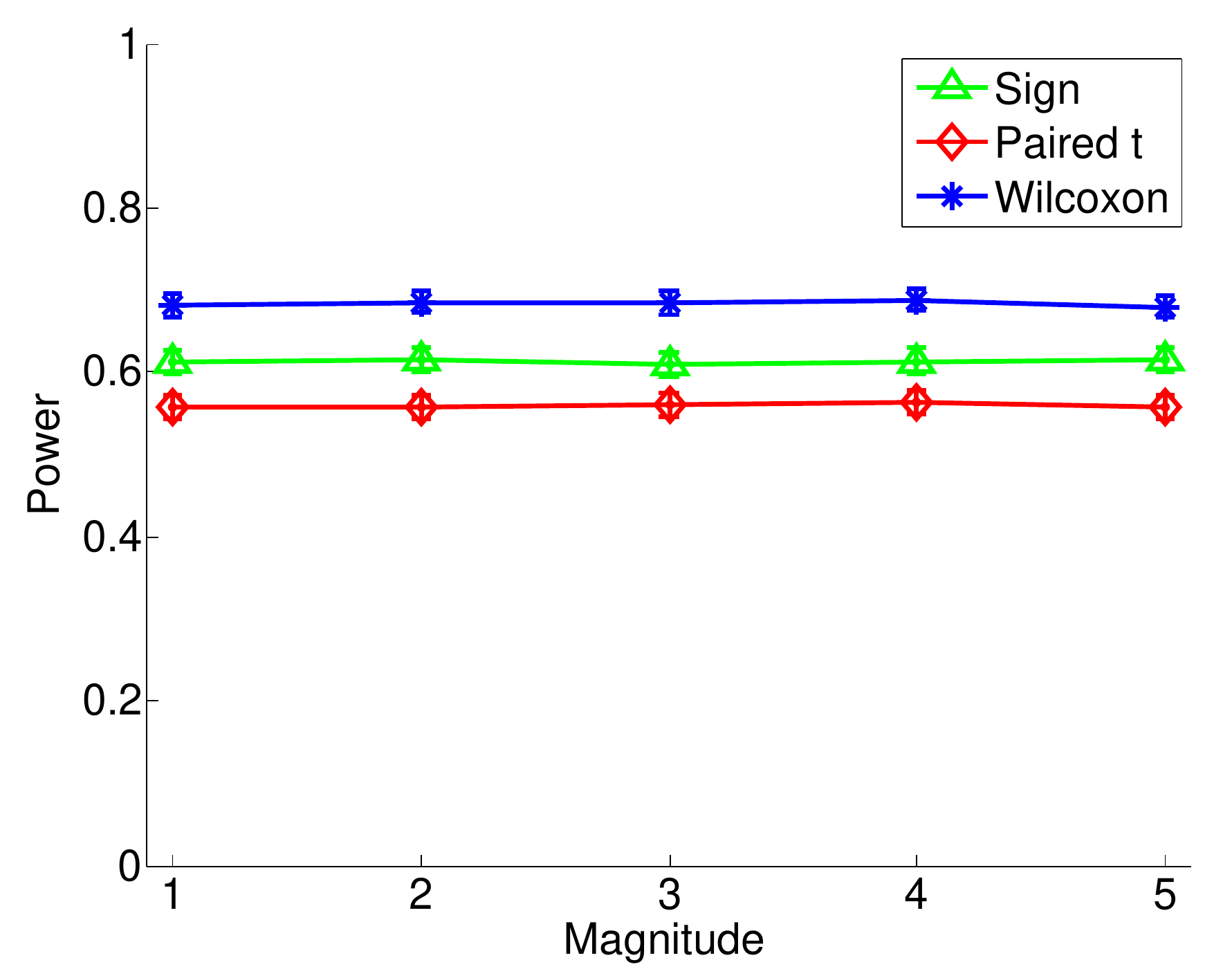}}\quad
\subfigure[Power versus~$c_v$.]{\includegraphics[width=0.3\textwidth]{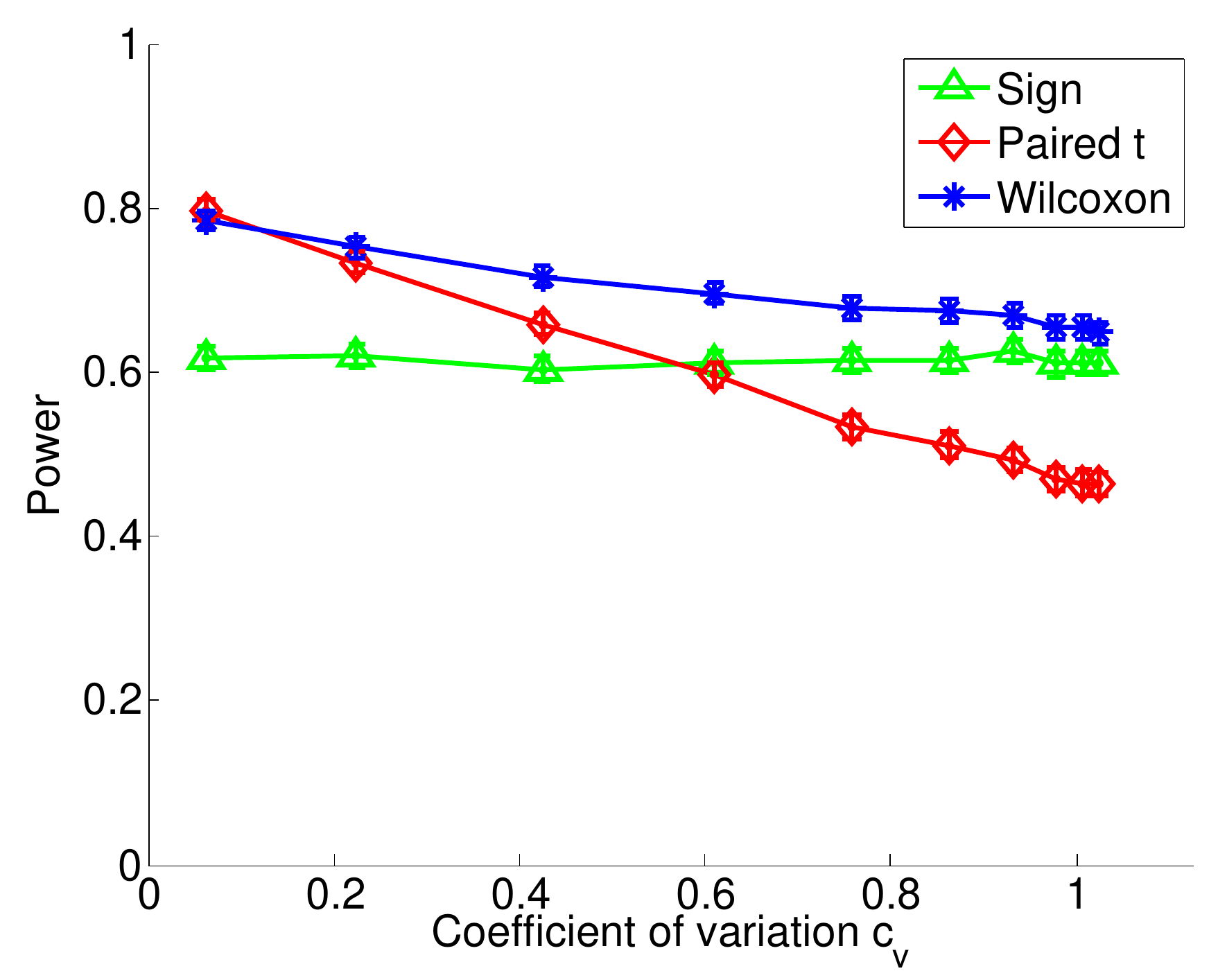}}\quad
\subfigure[Power versus~$c_v$.]{\includegraphics[width=0.3\textwidth]{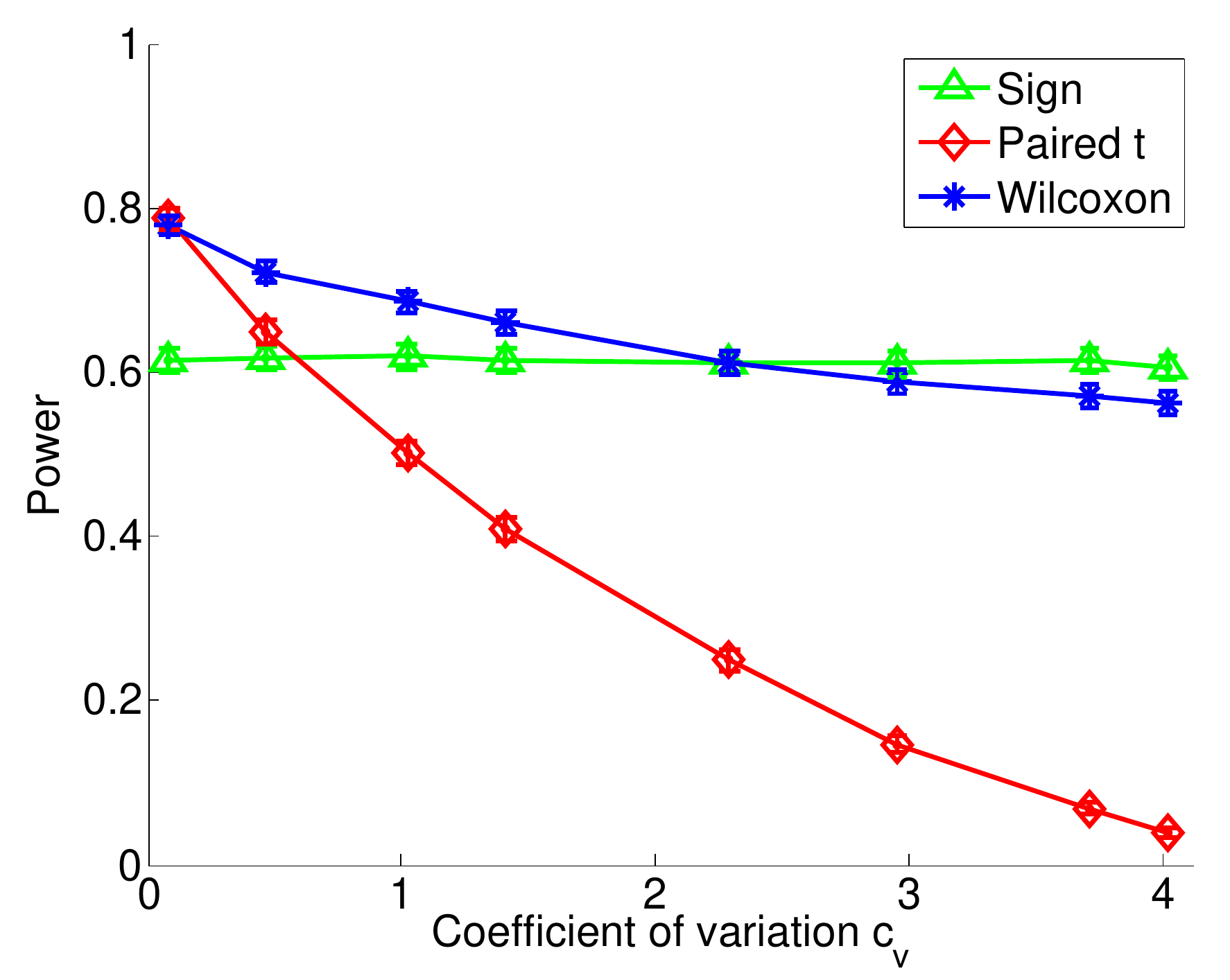}}
\caption{Effect of various parameters on the statistical power of the tests. \label{fig:Syn}}
\end{figure*}

\subsection{Comparison with Paired T-test}\label{sec:comp}
Besides the minimaxity of the sign test, it is interesting to identify (realistic) conditions on the nuisance parameters such that the sign test outperforms other popular tests. We demonstrate this by comparing the asymptotic power of the sign test with that of the paired t-test, whose test statistic is given in \eqref{eq:stats_pairt}.
Since both the sign test  and the  paired t-test only use  $\bY=\bX^B-\bX^A$ to compute the test statistics, we can only consider $\bY$ as the input, which is generated by $\mu_i$'s according to $Y_i \sim \mathcal{N}(\delta \mu_i, \mu_i^2)$, $\forall~i$. 
Let $m_1=\frac{1}{n}\sum_i \mu_i$ and $m_2=\frac{1}{n}\sum_i (\mu_i - m_1)^2$ be the mean and the variance of $\mu_i$'s. 
In the homogeneous case, $\mu_i$'s are the same and thus $m_2=0$. 
In the presence of within-group heterogeneity, however, $\mu_i$'s are different and $m_2$ may be large.  Hence, intuitively, it is reasonable to look at the \emph{coefficient of variation}~\cite{everitt2006cambridge}, $c_v=m_2 / m_1^2$, as a measure for within group heterogeneity.  
In fact, as  shown below, it is the determining factor for the testing performance.

\begin{thm}\label{thm:comp}
Let $n \rightarrow \infty$ and scale $\delta$ with $n$ such that $\delta \sqrt{n}$ remains constant\footnote{Such asymptotic scaling makes the power converge to some constant between $0$ and $1$, and thus making the power of the tests comparable \cite{van2000asymptotic}.}. Also assume that by increasing $n$ the values of $m_1$ and $m_2$ remain constant. Then, the asymptotic power of the $\alpha$-level two-sided sign test and the $\alpha$-level two-sided paired t-test are given by \eqref{eq:pw}, respectively. Moreover, the two-sided sign test has a larger asymptotic power if $c_v \geq \pi/2 -1$. 
\begin{align}\label{eq:pw}
&\textrm{Power of sign test: } Q\left( z_{\frac{\alpha}{2}} - \sqrt{\frac{2}{\pi}}\sqrt{n}\delta \right),\\&
\textrm{Power of paired t-test: }  Q\left( z_{\frac{\alpha}{2}} - \frac{\sqrt{n}\delta }{\sqrt{1+c_v}} \right),
\end{align}
where $Q(\cdot)$ is the tail function of the standard normal distribution and $z_{\frac{\alpha}{2}}$  is the $(1-\frac{\alpha}{2})$-th quantile of the standard normal distribution. 
\end{thm}
\textbf{Remark.} As shown in Theorem \ref{thm:comp}, $c_v$ is the key quantity that determines the performance of the paired t-test as compared to the signed test. The condition $c_v \geq \pi/2 -1$ is quite general, implying under a variety of the nuisance parameters, the sign test can outperform the paired t-test.

\section{NUMERICAL EXPERIMENTS} \label{sec:Exp}
In this section, we provide numerical evidence on the theoretical optimality and the practicality of our results. First, we evaluate our results on the synthetic data. Then, the viral infection disease dataset \cite{chen2012personal} is used to further evaluate the practicality of our theoretical findings. 

\subsection{Synthetic Data}
Here, we compare the performance of sign test with two other popular tests for paired data: the paired t-test and the Wilcoxon signed rank test \cite{van2000asymptotic}. 
All three tests calculate the test statistic using only $\bY$.
Hence we only consider different values of $\{\mu_i\}_{i=1}^n$ and generate the samples~$\bY$ according to $Y_i \sim \mathcal{N}(\delta \mu_i, \mu_i^2)$, for $i=1,...,n$.  
In all experiments, we fix the sample size $n=20$, $\delta= 3 / \sqrt{n}$, and the size of the tests $\alpha=0.05$.
We repeat experiments under each parameter setting $10,000$ times, and plot the corresponding $3~std$ confidence intervals.

Recall that, according to Theorem~\ref{thm:comp}, the coefficient of variation $c_v={m_2}/{m_1^2}$ determines the power of the paired t-test. 
Our first experiment examines if $c_v$ can quantify the within-group heterogeneity level reasonably for finite values of $n$. In this experiment, we generated $\bm{\mu}$ using the two-group model \cite{efron2008microarrays}, where the $\mu_i$'s are $50/50$ with values $1\times$ and $10\times$ of some given magnitudes.   
We plot the corresponding powers for different values of the given magnitudes while the corresponding $c_v$'s are kept fixed to the value $0.7$.
As shown in Fig. \ref{fig:Syn} (a), for all $3$ tests, the powers are the same for experiments under different magnitudes, implying that same value of $c_v$ always results in the same power regardless of the values of other parameters.  Hence, it is reasonable to assume that $c_v$ can well quantify the level of heterogeneity. 

In the other two experiments, we use $c_v$ to represent heterogeneity level and plot the power of different test versus different values of $c_v$. 
In Fig. \ref{fig:Syn} (b), the value of $\bm{\mu}$ is also generated by the two-group model as before. As can be seen in this figure, 
the sign test outperforms the paired t-test when $c_v$ exceeds $0.58$.  
This phenomenon is consistent with the condition in Theorem~\ref{thm:comp} where the threshold is computed as $c_v\geq \frac{\pi-2}{2}\approx0.57$. 
In Fig. \ref{fig:Syn} (c), $\bm{\mu}$ is generated according to the multi-group model with $5$ groups of different values.
As can be seen in the figure, the sign test has a better power than the Wilcoxon test when $c_v$ exceeds $2.3$. 
In addition, as $c_v$ increases, the power of both the paired t-test and the Wilcoxon signed test decreases, but the later decreases much slower than the former. 
This is in line with our intuition since the Wilcoxon signed-rank test statistic uses the sign information and is more robust to the within-group heterogeneity than the paired t-test statistic.  

\subsection{The Viral Infection Dataset}

\begin{table}
\centering
\begin{tabular}{|c|c|c|c|c|c|}
\hline
~ & Sign & Wil& Pair T & DESeq2 & Voom \\
\hline
Sign & 225 & 156 & 140 & 66 & 193\\
Wil &~&267&223&95&260\\
Pair T & ~ & ~ & 292 & 97 & 282  \\
DESeq2 & ~ & ~ & ~ & 170 & 163 \\ 
Voom & ~ & ~ & ~ & ~ & 628 \\  
\hline
\end{tabular}
\caption{Number of common discoveries across various methods: the sign test (Sign), Wilcoxon signed-rank test (Wil), paired t-test (Pair T), paired-mode DESeq2 (DESeq 2), paired-mode Voom (Voom). \label{tab:TestRes}}
\end{table}

\begin{figure}
\centering
\includegraphics[width=0.48\textwidth]{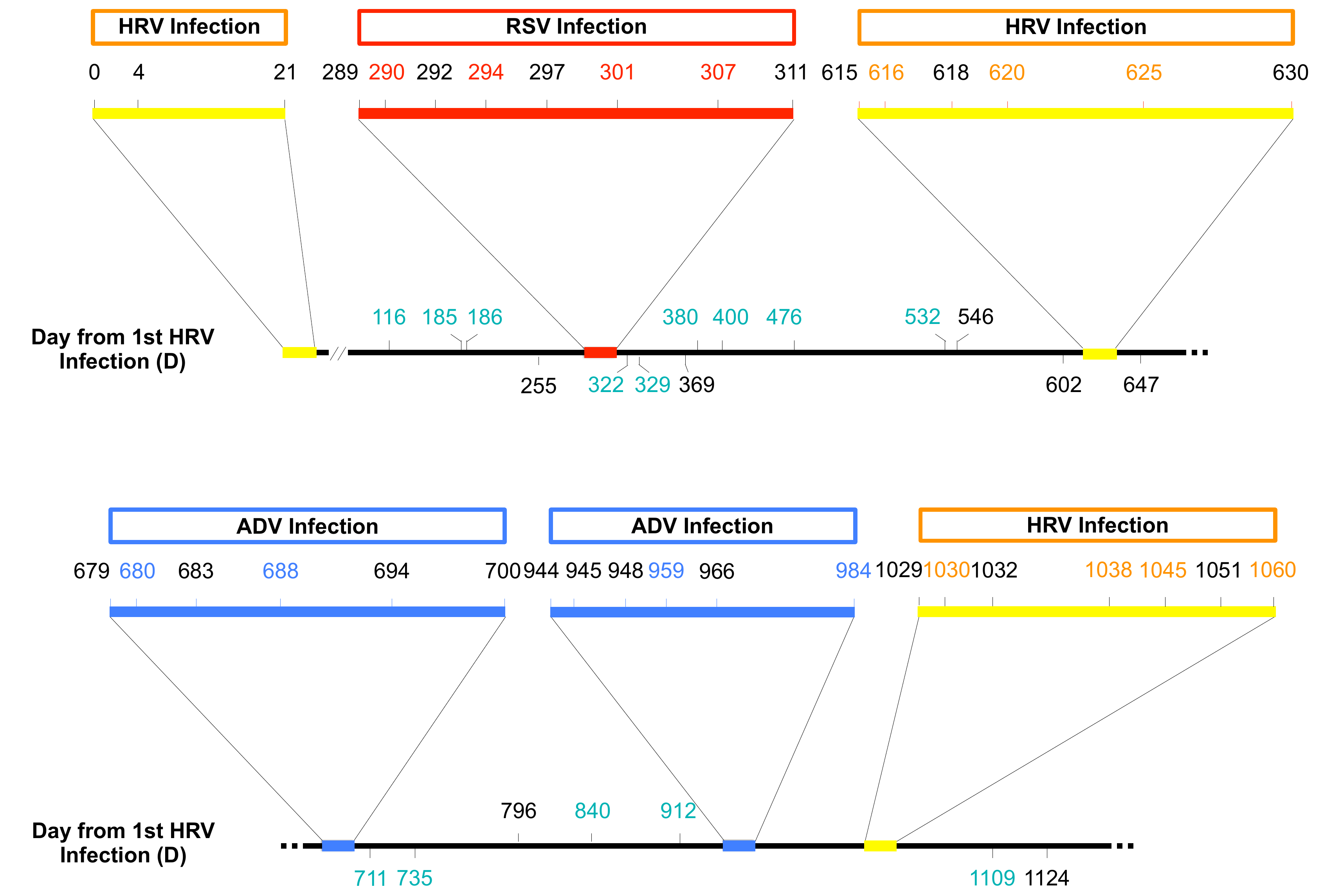}
\caption{The description of the VID dataset \citep{chen2012personal}. \label{fig:VirInf}}
\end{figure}

In VID \citep{chen2012personal}, one subject went through 6 viral infection periods in an overall time period of $1124$ days. 
During this period, $57$ RNA-Seq blood samples were collected under two conditions, healthy and sick;  see Fig. \ref{fig:VirInf} for more details. 
The task is to find differentially-expressed genes under the two conditions. 
We manually pair the samples under the two conditions that are close to each other in time; and altogether acquire $20$ data pairs. 
As shown in Fig. \ref{fig:EviDif} in the beginning of this manuscript, within-pair samples are more similar to each other than within-group samples. We compare the performance of the two-sided sign test \eqref{eq:SgnTst_twoside}, the Wilcoxon test \eqref{eq:stats_wil}, the paired t-test \eqref{eq:stats_pairt}.
We also report the performance of two popular differential expression analysis packages paired-mode DESeq2 \citep{love2014moderated} and paired-mode Voom \citep{law2014voom}.
According to \citep{pimentel2016differential}, they have the most promising performance among differential expression analysis tools. 
Prior to testing, genes with the total number of counts less than $50$ or having some counts less than or equal to $1$ are removed since we do not have enough observations of them.
For the sign test and the paired t-test, we used the size factor normalization method as in DESeq2. For all methods, after the p-value calculation, the BH procedure \citep{benjamini1995controlling} is used to control the false discovery rate (FDR) at the level of $0.1$.

\begin{figure}
\centering
\includegraphics[width=0.48\textwidth]{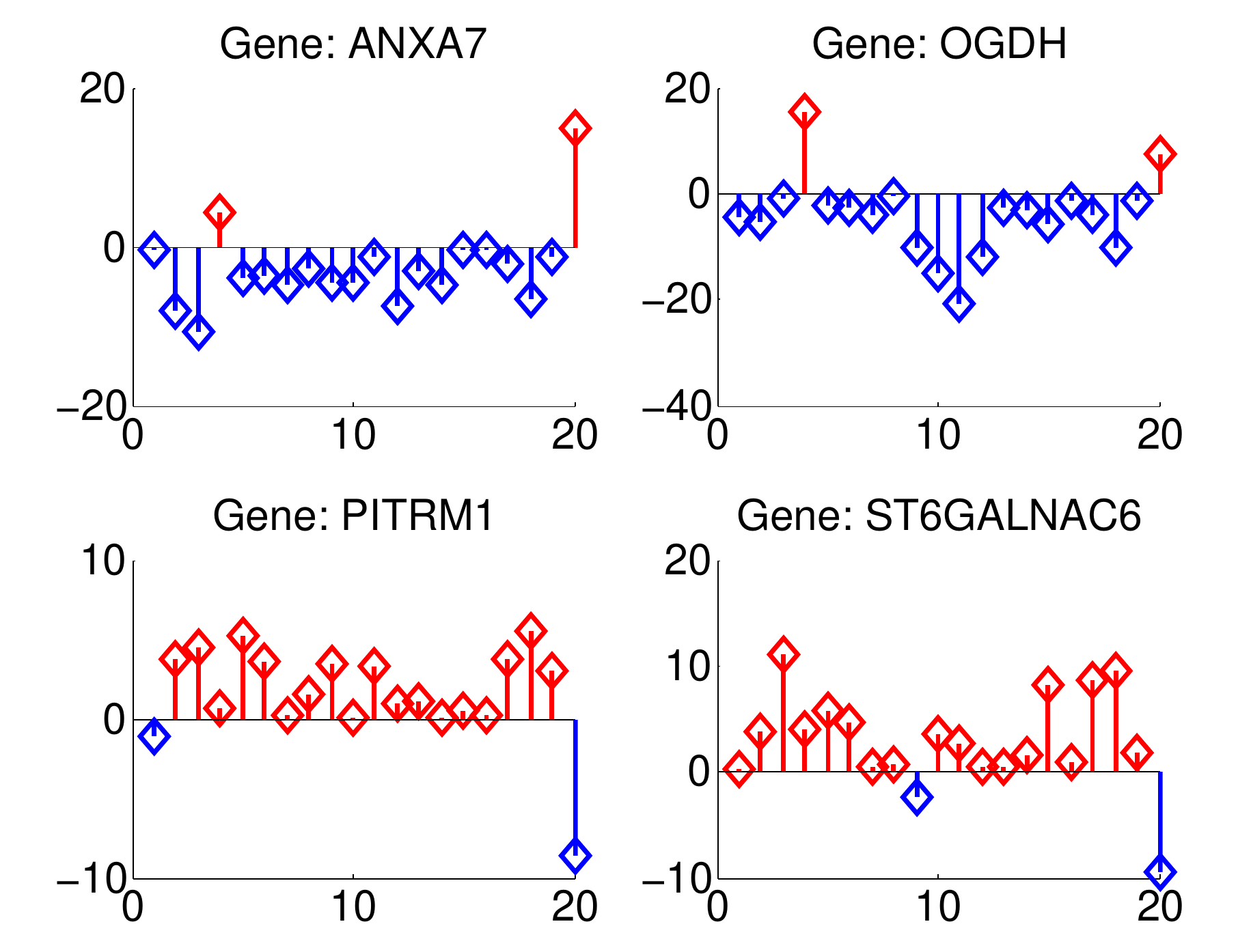}
\caption{Count difference of the top four genes discovered only by the sign test. X-axis: sample index; Y-axis: gene expression level. \label{fig:VirInf2}}
\end{figure}

We present our results in Table~\ref{tab:TestRes}, where the $ij$-th entry of this table is the number of genes discovered by both methods $i$ and $j$.   
There are $26$ genes discovered {\it only} by the sign test.
We plot the the paired differences $\bY=\bX^B-\bX^A$ for $4$ of these $26$ genes with smallest p-values in Fig. \ref{fig:VirInf2}, where the signals for differential expression are very strong.
For genes ANXA7, PITRM1, ST6GLLNAC6, sample $20$ has the opposite direction and a magnitude much larger than others. 
This prevents other methods that use the magnitude from discovering these genes.
But the sign test is robust to the heterogeneity in the magnitude and discovers them.

In Table~\ref{tab:TestRes}, the paired t-test has more discoveries than the sign test and the Wilcoxon test. A possible explanation is that the within-group heterogeneity level is not high enough due to the fact that the samples were all drawn from the same person. Voom makes more discoveries because it makes a strong assumption that a mean-variance function is shared across the genes (which also could lead to false discoveries).

\clearpage
\bibliographystyle{plain}
\bibliography{ref}


\clearpage
\onecolumn
\begin{center}
\textbf{\large Supplemental Materials}
\end{center}
\setcounter{section}{0}

\section{The details of the visualization \label{sec:vis}}

For any two samples $\mathbf{X}_i,\mathbf{X}_j$, we can calculate the histogram of the $k$ values $\log \vert \mathbf{X}_i-\mathbf{X}_j \vert  \in \mathbb{R}^k$. 
Concentration in small values indicates $\mathbf{X}_i,\mathbf{X}_j$ are similar. 
The histogram is averaged over all sample pairs for the within-pair difference and over all size-$2$ subsets of samples from the same group for the within-group difference.

\section{Proof of Fact \ref{fact:scl_inv} \label{sec:pf_fact}}
\begin{proof}
(Proof of Fact \ref{fact:scl_inv}) 
It suffices to show that 
\begin{equation}
\label{eq:Fact1}
\sup_{\gamma}\mathbb{E}_{0}[\phi]=\sup_{\gamma}\mathbb{E}_{0}[\psi],~~~~and ~~~~
\inf_{\gamma}\mathbb{E}_{1}[\phi]=\inf_{\gamma}\mathbb{E}_{1}[\psi]. 
\end{equation}
Now we prove the first equation in \eqref{eq:Fact1}. Let $\{Z_i^A, Z_i^B\}_{i=1}^n \overset{i.i.d.}{\sim}\mathcal{N}(0,1)$. Then by \eqref{eq:SampleModel} we can write 
\begin{equation*}
\mathbf{X}^{A}=\bm{\nu}+\bm{\mu}\circ\sqrt{\bm{\rho}}\circ\mathbf{Z}^{A},~
\mathbf{X}^{B}=\bm{\nu}+\delta\bm{\mu}+\bm{\mu}\circ\sqrt{1-\bm{\rho}}\circ\mathbf{Z}^{B}.
\end{equation*}
Similarly, we can write the transformed data as 
\begin{align*}
& \mathbf{a}+\mathbf{b}\circ\mathbf{X}^{A}=\bm{\nu}'+\bm{\mu}'\circ\sqrt{\bm{\rho}}\circ\mathbf{Z}{}^{A} \\
& \mathbf{a}+\mathbf{b}\circ\mathbf{X}^{B}=\bm{\nu}'+\delta\bm{\mu}'+\bm{\mu}'\circ\sqrt{1-\bm{\rho}}\circ\mathbf{Z}^{B},
\end{align*}
for $\bm{\nu}' \triangleq \mathbf{a}+\mathbf{b}\circ\bm{\nu}$ and $\bm{\mu}' \triangleq \mathbf{b}\circ\bm{\mu}$. 

The key idea for proving the first equation in \eqref{eq:Fact1} is that $\{\bm{\nu}, \bm{\mu}\}$ and $\{\bm{\nu}',\bm{\mu}'\}$ actually consist of the same parameter space $\mathbb{R}^n\times \mathbb{R}_{>0}^n$. 
To be more exact, notice that under the null hypothesis $\delta=0$. We have  
\begin{align*}
 & \sup_{\bm{\nu}, \bm{\mu},\bm{\rho}}\E_{0}[\psi(\mathbf{\mathbf{X}}^{A},\mathbf{X}^{B})]\\
  =&\sup_{\bm{\nu}, \bm{\mu},\bm{\rho}}\mathbb{E}[\phi(\bm{\nu}'+\bm{\mu}'\circ\sqrt{\bm{\rho}}\circ\mathbf{Z}{}^{A},\bm{\nu}'+\bm{\mu}'\circ\sqrt{1-\bm{\rho}}\circ\mathbf{Z}^{B})]\\
 =&\sup_{\bm{\nu}',\bm{\mu}',\bm{\rho}}\mathbb{E}[\phi(\bm{\nu}'+\bm{\mu}'\circ\sqrt{\bm{\rho}}\circ\mathbf{Z}{}^{A},\bm{\nu}'+\bm{\mu}'\circ\sqrt{1-\bm{\rho}}\circ\mathbf{Z}^{B})]\\
 =&\sup_{\bm{\nu}, \bm{\mu},\bm{\rho}}\mathbb{E}[\phi(\bm{\nu}+\bm{\mu}\circ\sqrt{\bm{\rho}}\circ\mathbf{Z}{}^{A},\bm{\nu}+\bm{\mu}\circ\sqrt{1-\bm{\rho}}\circ\mathbf{Z}^{B})]\\
 =&\sup_{\bm{\nu}, \bm{\mu},\bm{\rho}}\mathbb{E}_{0}[\phi(\mathbf{\mathbf{X}}^{A},\mathbf{X}^{B})],
\end{align*}
where in the second equation we use the fact that $\{\bm{\nu}, \bm{\mu}\}$ and $\{\bm{\nu}',\bm{\mu}'\}$ consist of the same parameter space. The second equation of \eqref{eq:Fact1} can be shown similarly. Then $\phi\doteq\psi$. 
\end{proof}

\section{Proof of Theorems}
\subsection{Proof of Theorem \ref{thrm:opt_two_side} \label{pf:opt_two_side}}
\begin{proof} (Proof of Theorem \ref{thrm:opt_two_side})
First, for this problem, we should restrict ourselves to symmetrical tests, i.e. any $\phi$ such that $\phi(\bX^A,\bX^B)=\phi(-\bX^A,-\bX^B)$. This is because for any $(\bX^A,\bX^B)$, the distribution $(-\bX^A,-\bX^B)$ is also valid for the maximin problem. 
Second, according to Lemma \ref{lem:two_to_one}, it suffices to consider the symmetrical tests using only $\bY$.  

Now consider any $\alpha$-level symmetrical test that uses only $\bY$, namely $\phi(\bY)$. 
For any set of nuisance parameters $\gamma$, let $f_{0, \gamma}(\cdot)$, $f_{1, \gamma}(\cdot)$ 
be the density function under the null and the alternative hypothesis respectively. 
Define 
\begin{align*}
\mathcal{Y}^+ = \{\mathbf{y}: f_{1, \gamma} (\mathbf{y}) > f_{1, \gamma} (\mathbf{-y}) \}, \\
\mathcal{Y}^- = \{\mathbf{y}: f_{1, \gamma} (\mathbf{y}) < f_{1, \gamma} (\mathbf{-y}) \}. 
\end{align*}
Then by explicitly writing out the density function, it is not hard to see that for any $\mathbf{y}\in \mathcal{Y}^+$, we have $-\mathbf{y}\in \mathcal{Y}^-$. As the null distribution is symmetrical around $0$, we have
\begin{equation} \label{eq:two_sided size}
\sup_\gamma \E_0[\phi(\bY) \ind_{\{\bY\in \mathcal{Y}^+\}}] = \sup_\gamma  \E_0[\phi(\bY) \ind_{\{\bY\in \mathcal{Y}^-\}}] \leq \frac{\alpha}{2}
\end{equation}

Next consider the power of $\phi$. As $\phi(\bY) \ind_{\{\bY\in \mathcal{Y}^+\}}$ is a $\frac{\alpha}{2}$-level test, by the optimality of the one-sided sign test, $\forall~\epsilon>0$, $\exists~\gamma^*$, s.t.,
\begin{align}\label{eq:two_sided_eq1}
\E_{1,\gamma^*} [\phi(\bY) \ind_{\{\bY\in \mathcal{Y}^+\}}] -\epsilon  < \inf_\gamma \E_1[\phi(\bY) \ind_{\{\bY\in \mathcal{Y}^+\}}] \leq \inf_{\gamma} \E_1 [\tilde{\phi}^S(\bY)]. 
\end{align}
Also we have
\begin{align}\label{eq:two_sided_eq2}
&\E_{1,\gamma^*} [\phi(\bY) \ind_{\{\bY\in \mathcal{Y}^-\}}] = \int_{\mathcal{Y}^-} \phi(\mathbf{y}) f_{1, \gamma^*} (\mathbf{y}) d\mathbf{y}
= \int_{\mathcal{Y}^-} \phi(\mathbf{y}) f_{0, \gamma^*}(\mathbf{y}) \frac{f_{1, \gamma^*} (\mathbf{y})}{f_{0, \gamma^*} (\mathbf{y})} d\mathbf{y} \nonumber\\
& \leq \sup_{\mathbf{y}\in \mathcal{Y}^-} \left( \frac{f_{1,\gamma^*} (\mathbf{y})}{f_{0,\gamma^*} (\mathbf{y})} \right)\sup_{\gamma}  \E_0[\phi(\bY) \ind_{\{\bY\in \mathcal{Y}^-\}}]
\leq \frac{\alpha}{2}  \exp(-\frac{\delta^2 n}{2}). 
\end{align}
In the last inequality of \eqref{eq:two_sided_eq2}, the first term is due to the fact that if $\mathbf{y}\in \mathcal{Y}^-$, we have $\frac{f_{1, \gamma^*} (\mathbf{y})}{f_{1, \gamma^*} (-\mathbf{y})}\leq 1$, giving $\delta \sum_{i=1}^i \frac{y_i}{\mu_i}\leq 0$, which further gives $\forall~\mathbf{y}\in \mathcal{Y}^-$,
\begin{equation*}
\frac{f_{1, \gamma^*} (\mathbf{y})}{f_{0, \gamma^*} (\mathbf{y})} = \exp(-\frac{n\delta^2}{2}+ \delta \sum_{i=1}^i \frac{y_i}{\mu_i}) \leq \exp(-\frac{n\delta^2}{2}). 
\end{equation*}
The second term is due to (\ref{eq:two_sided size}). 
Finally, combining \eqref{eq:two_sided_eq1} and \eqref{eq:two_sided_eq2}, we reach that the power of $\phi$ 
\begin{align*}
&\inf_{\gamma} \E_1 [\phi(\bY)]  \leq \E_{1,\gamma^*}[\phi(\bY)] = \E_{1,\gamma^*} [\phi(\bY) \ind_{\{\bY\in \mathcal{Y}^+\}} + \phi(\bY) \ind_{\{\bY\in \mathcal{Y}^-\}}]\\
& \leq \inf_{\gamma} \E_1 [\tilde{\phi}^S(\bY)] + \frac{\alpha}{2}  \exp(-\frac{n \delta^2}{2}) +\epsilon \leq  \inf_{\gamma} \E_1 [\phi^S(\bY)] + \frac{\alpha}{2}  \exp(-\frac{n \delta^2}{2}) + \epsilon, 
\end{align*}
where we recall that $\phi^S$ is the two-sided sign test as defined in \eqref{def:two_side_sign}.
Finally, let $\epsilon\to 0$ to complete the proof.
\end{proof}

\subsection{Proof of Theorem \ref{thm:comp}}
\begin{proof}
(Proof of Theorem \ref{thm:comp})
Let us start by defining the following notations: we will use $\covind$ and $\covinp$ to denote convergence in distribution and convergence in probability, respectively.
Next, we first compute the power of the two-sided sign test. Let us consider the sufficient statistics  $W=\sum_i \ind_{\{Y_i>0\}}$ for the sign test. Clearly, under the null distribution, $\frac{W-0.5 n}{\sqrt{0.25n}} \covind \mathcal{N} (0,1)$ due to the standard central limit theorem. Thus the two-sided sign test asymptotically rejects when $\vert \frac{W-0.5 n}{\sqrt{0.25n}} \vert \geq z_{\frac{\alpha}{2}}$. Furthermore, under the alternative hypothesis, by Taylor's expansion, we have
\begin{equation*}
\theta=\Prob(Y_i \geq 0) = Q(-\delta) = 0.5 +\frac{1}{\sqrt{2\pi}} \delta + O(\delta^2).  
\end{equation*}
Without loss of generality assume that $\delta>0$. Then the power of the test can be computed as 
\[
 \Prob (\vert \frac{W-0.5 n}{\sqrt{0.25n}} \vert \geq z_{\frac{\alpha}{2}}) \approx  \Prob ( \frac{W-0.5 n}{\sqrt{0.25n}} \geq z_{\frac{\alpha}{2}}),  \\
\]
where by using the approximately equality we neglect the lower tail, a small quantity that decreases exponentially with $\alpha$. This power can be further computed as 
\begin{equation} \label{eq:pfpt1}
\begin{split}
&  Q\left( z_{\frac{\alpha}{2}} - \sqrt{\frac{2}{\pi}}\sqrt{n}\delta \right).
\end{split}
\end{equation}

On the other hand, for the paired t-test, we can asymptotically write the numerator and the denominator of the $T$-statistics as
\begin{align*}
& \sqrt{n} \bar{Y} \covind \mathcal{N} \left(\sqrt{n} \delta m_1, m_1^2 + m_2 \right)\\
& \frac{1}{n-1} \sum_i (Y_i - \bar{Y})^2 \covinp m_1^2 + m_2,
\end{align*}
where we recall that $\sqrt{n} \delta$ is some given constant by appropriate scaling of $\delta$.  Therefore, under the null distribution, we have $T \covind \mathcal{N}(0,1)$. Consequently, the paired t-test rejects when $\vert T \vert \geq z_{\frac{\alpha}{2}}$. Again by assuming $\delta>0$ under the alternative hypothesis, the power of the test can be written as
\small
\begin{equation}\label{eq:pfpt2}
\begin{split}
\Prob(\vert T \vert \geq z_{\frac{\alpha}{2}}) \approx \Prob( T \geq z_{\frac{\alpha}{2}})  
\to Q\left( z_{\frac{\alpha}{2}} - \frac{\sqrt{n}\delta }{\sqrt{1+c_v}}\right). 
\end{split}
\end{equation}
\normalsize
Combining \eqref{eq:pfpt1} and \eqref{eq:pfpt2} will complete the proof. 
\end{proof}

\section{Proofs of Lemmas}
\subsection{Proof of Lemma \ref{lem:two_to_one}}

\begin{proof}
Since $\mathbf{X}^{B}=\mathbf{X}^{A}+\mathbf{Y}$, $\phi$ can be equivalently represented as $\tilde{\phi}(\bY,\bX^A) = \phi (\bX^A, \bY+\bX^A)$. 
Let $\psi(\mathbf{Y})=\tilde{\phi}(\mathbf{Y},0)$.
If $\phi$ is symmetric, then 
\begin{equation}
\begin{split}
& \psi(\mathbf{Y}) = \tilde{\phi}(\mathbf{Y},0) = \phi (0, \bX^B-\bX^A) \\
& = \phi (0, \bX^A-\bX^B) = \tilde{\phi}(-\mathbf{Y},0) = \psi(-\mathbf{Y}),
\end{split}
\end{equation}
giving that $\psi$ is also symmetric. 

We next show that $\tilde{\phi}\dotleq\psi$. 
First notice that for each $i$, $(X_{i}^{A},Y_{i})$ follows a joint Gaussian distribution:
\[
\left[\begin{array}{c}
X_{i}^{A}\\
Y_{i}
\end{array}\right]
\sim{\cal N}\left(\begin{array}{cc}
\left[\begin{array}{c}
\nu_{i}\\
\delta\mu_{i}
\end{array}\right], & \left[\begin{array}{cc}
\rho_{i}\mu_{i}^{2} & -\rho_i\mu_{i}^{2}\\
-\rho_i\mu_{i}^{2} & \mu_{i}^{2}
\end{array}\right]\end{array}\right),
\]
and the samples are independent across all indices $i=1,\ldots,n$. 
Then $X_{i}^{A}\vert Y_i=y_{i}\sim{\cal N}(\nu_{i}-\rho_{i}y_{i}+\rho_{i}\delta\mu_{i},\rho_{i}(1-\rho_{i})\mu_{i}^{2})$ and $\forall~i\neq j$, $X_{i}^{A}\bot X_{j}^{A}\vert\mathbf{Y}$, where we note that $\left[X_{i}^{A}\vert Y_i=y_{i}, \nu_i=\rho_i=0\right] =0 $.

For the worst-case size, 
\begin{equation}
\begin{split}\label{eq:lemmaXtoYtemp1}
&\sup_{\bm{\mu,\bm{\nu},\bm{\rho}}}\E_{0}[\phi]  =\sup_{\bm{\mu,\bm{\nu},\bm{\rho}}}\E_{0,\mathbf{Y}}[\E_{0,\mathbf{X}^{A}\vert\mathbf{Y}}[\tilde{\phi}(\mathbf{Y},\mathbf{X}^{A})\vert\mathbf{Y}]] \\& \geq \sup_{\bm{\mu,\nu=\rho=0}}\E_{0,\mathbf{Y}}[\E_{0,\mathbf{X}^{A}\vert\mathbf{Y}}[\tilde{\phi}(\mathbf{Y},\mathbf{X}^{A})\vert\mathbf{Y}]] \\&= \sup_{\bm{\mu,\nu=\rho=0}}\E_{0,\mathbf{Y}}[\E_{0,\mathbf{X}^{A}\vert\mathbf{Y}}[\tilde{\phi}(\mathbf{Y},0)\vert\mathbf{Y}]] \\&= \sup_{\bm{\mu,\nu=\rho=0}}\E_{0,\mathbf{Y}}[\E_{0,\mathbf{X}^{A}\vert\mathbf{Y}}[\psi(\mathbf{Y})\vert\mathbf{Y}]] = \sup_{\bm{\mu,\nu,\rho}}\E_0[\psi],
\end{split}
\end{equation}
where the last equality is because $\psi(\mathbf{Y})$ does not depend on $\bm{\nu},\bm{\rho}$. Similarly, the power of the two tests can be related by

\begin{equation}
\begin{split}  \label{eq:lemmaXtoYtemp2}
& \inf_{\bm{\mu,\bm{\nu},\bm{\rho}},s_\delta}\E_{1}[\phi] = \inf_{\bm{\mu,\bm{\nu},\bm{\rho}},s_\delta}\E_{1,\mathbf{Y}}[\E_{1,\mathbf{X}^{A}\vert\mathbf{Y}}[\tilde{\phi}(\mathbf{Y},\mathbf{X}^{A})\vert\mathbf{Y}]] \\&
\leq \inf_{\bm{\mu,\nu=\rho=0},s_\delta}\E_{1,\mathbf{Y}}[\E_{1,\mathbf{X}^{A}\vert\mathbf{Y}}[\tilde{\phi}(\mathbf{Y},\mathbf{X}^{A})\vert\mathbf{Y}]] \\&
= \inf_{\bm{\mu,\nu=\rho=0},s_\delta}\E_{1,\mathbf{Y}}[\E_{1,\mathbf{X}^{A}\vert\mathbf{Y}}[\tilde{\phi}(\mathbf{Y},0)\vert\mathbf{Y}]] \\&
= \inf_{\bm{\mu,\nu=\rho=0},s_\delta}\E_{1,\mathbf{Y}}[\E_{1,\mathbf{X}^{A}\vert\mathbf{Y}}[\psi(\mathbf{Y})\vert\mathbf{Y}]] = \inf_{\bm{\mu,\nu, \rho},s_\delta}\E_1 [\psi],
\end{split}
\end{equation}
where $s_\delta=1$ for the one-sided case and $s_\delta\in\{1,-1\}$ for the two-sided case. Combining \eqref{eq:lemmaXtoYtemp1} and \eqref{eq:lemmaXtoYtemp2} we conclude $\phi\dotleq\psi$. 
\end{proof}

\subsection{Proof of Lemma \ref{lem:stepopt}}
\begin{proof}
Let $Z_1,\cdots,Z_n \overset{i.i.d.}{\sim} \mathcal{N}(0,1)$, and let $P_0$, $P_1$ denote the $n$-D product probability measure of $\mathbf{Z}=(Z_1,\cdots,Z_n)$ and $\mathbf{Z}+\delta=(Z_1+\delta,\cdots,Z_n+\delta)$ respectively. 
In addition, denote the probability of orthant $\mathbf{o}$ under measure $P_0$ and $P_1$ by $P_0^{\mathbf{o}}$ and $P_1^{\mathbf{o}}$ respectively. 
Notice that $\sum_{\mathbf{o}\in\mathcal{O}} P_0^{\mathbf{o}}=\sum_{\mathbf{o}\in\mathcal{O}} P_1^{\mathbf{o}}=1$, and for any $\mathbf{o}\in\mathcal{O}$, $P_0^{\mathbf{o}}=2^{-n}$. 

Depending only on the sign of the data, the sign test has constant value on each orthant. 
In fact, it is not hard to see that the level-$\alpha$ one-sided sign test $\phi^S(\mathbf{Y})$ \eqref{eq:SgnTst_oneside} maximizes its power by assigning $1$ to orthants with larger values of $P_1^{\mathbf{o}}$ until the corresponding size reaches $\alpha$. 
We next show that such procedure has a power that upper bounds the power of any level-$\alpha$ test in $\mathcal{S}$, which proves the lemma. 

\textbf{We first define some useful quantities.}
Now consider any $\phi \in \mathcal{S}(\omega)$. 
Define the discretized space of $\bm{\mu}$ to be $D=\{\bm{\mu}=(\mu_1,\cdots,\mu_n): \forall i, \mu_i=(1+\omega)^{d_i}, d_i\in\mathbb{Z}\}$. 
Consider any $\bm{\mu}\in D$ with $\mathbf{d}=(d_1,\cdots,d_n)$.
For any $n$-D box $I^\mathbf{o}_{\mathbf{b}}$, the element-wise multiplication by $\bm{\mu}$ maps it to $I^\mathbf{o}_{\mathbf{b}+\mathbf{d}}$.
Define $f(\mathbf{z})$ to be the density function of $\mathbf{Z}$.
Let $\alpha(\bm{\mu})=\E_{0}[\phi(\mathbf{Y})]$ and $\beta(\bm{\mu})=\E_{1}[\phi(\mathbf{Y})]$ be the size and power of $\phi$ when the nuisance parameters have value $\bm{\mu}$.
We have
\begin{align*}
& \alpha(\bm{\mu})  =\E[\phi(\bm{\mu} \circ \mathbf{Z})]=\int_{\mathbb{R}^n}\phi(\bm{\mu} \circ \mathbf{z})f(\mathbf{z})d\mathbf{z} = \sum_{\mathbf{o}\in\mathcal{O}}\sum_{-\infty<\mathbf{b}<\infty}\phi_{\mathbf{d}+\mathbf{b}}^{\mathbf{o}}P_{0}(I_{\mathbf{b}}^{\mathbf{o}}), \\
& \beta(\bm{\mu}) =\E[\phi(\bm{\mu}\circ(\mathbf{Z}+\delta))]=\int_{\mathbb{R}^n}\phi(\bm{\mu} \circ \mathbf{z})f(\mathbf{z}-\delta)d\mathbf{z} =\sum_{\mathbf{o}\in\mathcal{O}} \sum_{-\infty<\mathbf{b}<\infty}\phi_{\mathbf{d}+\mathbf{b}}^{\mathbf{o}}P_{1}(I_{\mathbf{b}}^{\mathbf{o}}).
\end{align*}
Moreover, we can write the size and power corresponding to orthant $\mathbf{o}$ as 
\begin{align} \label{eq:stepopt_pf_1}
\alpha^\mathbf{o}(\bm{\mu}) = \sum_{-\infty<\mathbf{b}<\infty}\phi_{\mathbf{d}+\mathbf{b}}^{\mathbf{o}}P_{0}(I_{\mathbf{b}}^{\mathbf{o}}),~~~~~~~~
\beta^\mathbf{o}(\bm{\mu}) = \sum_{-\infty<\mathbf{b}<\infty}\phi_{\mathbf{d}+\mathbf{b}}^{\mathbf{o}}P_{1}(I_{\mathbf{b}}^{\mathbf{o}}),
\end{align}
where we note that $\alpha(\bm{\mu}) =\sum_{\mathbf{o}\in\mathcal{O}} \alpha^\mathbf{o}(\bm{\mu})$ and $\beta(\bm{\mu}) =\sum_{\mathbf{o}\in\mathcal{O}} \beta^\mathbf{o}(\bm{\mu})$. Next, for each orthant $\mathbf{o}$ and any positive integer $m$, define the $m$-th approximation of the probability measure $P_0^{\mathbf{o}}$, $P_1^{\mathbf{o}}$, the size $\alpha^\mathbf{o}(\bm{\mu})$ and the power $\beta^\mathbf{o}(\bm{\mu})$ to be  
\begin{align}\label{eq:stepopt_pf_2}
P_{0,m}^{\mathbf{o}} &=\sum_{-m\leq\mathbf{b}\leq m} P_0(I^\mathbf{o}_\mathbf{b}),&P_{1,m}^{\mathbf{o}} &= \sum_{-m\leq\mathbf{b}\leq m} P_1(I^\mathbf{o}_\mathbf{b})\\ \label{eq:stepopt_pf_3}
\alpha_{m}^\mathbf{o}(\bm{\mu}) &= \sum_{-m\leq\mathbf{b}\leq m}\phi_{\mathbf{d}+\mathbf{b}}^{\mathbf{o}}P_{0}(I_{\mathbf{b}}^{\mathbf{o}}), &\beta^\mathbf{o}_m(\bm{\mu}) &= \sum_{-m\leq\mathbf{b}\leq m}\phi_{\mathbf{d}+\mathbf{b}}^{\mathbf{o}}P_{1}(I_{\mathbf{b}}^{\mathbf{o}}),
\end{align}
where the summation is over all indices $\mathbf{b}\in \mathbb{Z}^n$ with elements all between $-m$ and $m$. 
It is not hard to see that $P_{0,m}^{\mathbf{o}} \uparrow P_{0}^{\mathbf{o}}$, $P_{1,m}^{\mathbf{o}} \uparrow P_{1}^{\mathbf{o}}$, $\alpha_{m}^\mathbf{o}(\bm{\mu})\uparrow\alpha^\mathbf{o}(\bm{\mu})$, and $\beta_{m}^\mathbf{o}(\bm{\mu})\uparrow\beta^\mathbf{o}(\bm{\mu})$. 

\textbf{We next show the key step of the proof:} for every $m$, 
\begin{align} \label{eq:stepopt_pf_4}
\inf_{\bm{\mu}\in D} \sum_{\mathbf{o}\in\mathcal{O}} \frac{\beta_{m}^\mathbf{o}(\bm{\mu})}{P_{1,m}^{\mathbf{o}} } \leq \sup_{\bm{\mu}\in D} \sum_{\mathbf{o}\in\mathcal{O}} \frac{\alpha_{m}^\mathbf{o}(\bm{\mu})}{P_{0,m}^{\mathbf{o}}}.
\end{align}

Recall that $\bm{\mu}=(\mu_1,\cdots,\mu_n)$, for $\mu_i=(1+\omega)^{d_i}$. 
We can write,
\begin{align*}
&\sup_{\bm{\mu}\in D} \sum_{\mathbf{o}\in\mathcal{O}} \frac{\alpha_{m}^\mathbf{o}(\bm{\mu})}{P_{0,m}^{\mathbf{o}}} - \inf_{\bm{\mu}\in D} \sum_{\mathbf{o}\in\mathcal{O}} \frac{\beta_{m}^\mathbf{o}(\bm{\mu})}{P_{1,m}^{\mathbf{o}} } =\lim_{l \to \infty } \left[ \sup_{-l\leq\mathbf{d}\leq l} \sum_{\mathbf{o}\in\mathcal{O}} \frac{\alpha_{m}^\mathbf{o}(\bm{\mu})}{P_{0,m}^{\mathbf{o}}} - \inf_{-l\leq\mathbf{d}\leq l} \sum_{\mathbf{o}\in\mathcal{O}} \frac{\beta_{m}^\mathbf{o}(\bm{\mu})}{P_{1,m}^{\mathbf{o}} }\right]\\
& = \lim_{l \to \infty } \left[\sup_{-l\leq\mathbf{d}\leq l} \sum_{\mathbf{o}\in\mathcal{O}} \frac{\sum_{-m\leq\mathbf{b}\leq m}\phi_{\mathbf{d}+\mathbf{b}}^{\mathbf{o}}P_{0}(I_{\mathbf{b}}^{\mathbf{o}})}{P_{0,m}^{\mathbf{o}}} - \inf_{-l\leq\mathbf{d}\leq l} \sum_{\mathbf{o}\in\mathcal{O}} \frac{\sum_{-m\leq\mathbf{b}\leq m}\phi_{\mathbf{d}+\mathbf{b}}^{\mathbf{o}}P_{1}(I_{\mathbf{b}}^{\mathbf{o}})}{P_{1,m}^{\mathbf{o}} } \right]\\
& \geq \lim_{l \to \infty }\frac{1}{(2l+1)^n} \left[ \sum_{-l\leq \mathbf{d} \leq l} \sum_{\mathbf{o}\in\mathcal{O}} \frac{\sum_{-m\leq\mathbf{b}\leq m}\phi_{\mathbf{d}+\mathbf{b}}^{\mathbf{o}}P_{0}(I_{\mathbf{b}}^{\mathbf{o}})}{P_{0,m}^{\mathbf{o}}} -  \sum_{-l\leq \mathbf{d} \leq l} \sum_{\mathbf{o}\in\mathcal{O}} \frac{\sum_{-m\leq\mathbf{b}\leq m}\phi_{\mathbf{d}+\mathbf{b}}^{\mathbf{o}}P_{1}(I_{\mathbf{b}}^{\mathbf{o}})}{P_{1,m}^{\mathbf{o}} } \right]\\
&\geq \lim_{l \to \infty } O(\frac{m}{l}) =0,
\end{align*}
which gives \eqref{eq:stepopt_pf_4}. 
The first equality is because $D = \lim_{l\to\infty} \{ \bm{\mu}: \mu_i=(1+\omega)^{d_i}, -l\leq \mathbf{d}\leq l\}$. 
The second equality is because of \eqref{eq:stepopt_pf_3}.
The first inequality is obtained by replacing $\sup$ and $\inf$ by averaging.
The second inequality is a little tricky.
Inside the square brackets, both the first and the second big term can be rearranged by $\phi^{\mathbf{o}}_{\mathbf{j}}$ to have the form $\sum_{-l-m\leq \mathbf{j} \leq l+m} c_{\mathbf{j}} \phi^{\mathbf{o}}_{\mathbf{j}}$.
A careful inspection reveals that $c_{\mathbf{j}}=1$ for any $-l+m\leq \mathbf{j} \leq l-m$. 
Indeed, e.g. for the first big term, as long as $-l+m\leq \mathbf{j} \leq l-m$, $\phi^{\mathbf{o}}_{\mathbf{j}}$ will be multiplied by $P_{0}(I_{\mathbf{b}}^{\mathbf{o}})$ once for each $-m\leq \mathbf{b} \leq m$, and the sum of those coefficients, $\sum_{-m\leq \mathbf{b} \leq m} P_{0}(I_{\mathbf{b}}^{\mathbf{o}})$, exactly equals the denominator $P_{0,m}^{\mathbf{o}}$ according to \eqref{eq:stepopt_pf_2}. 
Thus, for all $-l+m\leq \mathbf{j} \leq l-m$, $\phi^{\mathbf{o}}_{\mathbf{j}}$ will have the same coefficients for both the first big term and the second big term, resulting altogether $(2l-2m+1)^n$ terms canceling each other.
As a result, there remain $(2l+2m+1)^n-(2l-2m+1)^n = O(m l^{n-1})$ terms, whose corresponding summation can be upper bounded by $O(\frac{m}{l})$ because of the multiplicative factor $\frac{1}{(2l+1)^n}$ outside the square brackets.
Finally, $O(\frac{m}{l})$ vanishes as $l\to\infty$.

\textbf{Next we prove the limiting case of \eqref{eq:stepopt_pf_4}:}  
\begin{align} \label{eq:stepopt_pf_5}
\inf_{\bm{\mu}\in D} \sum_{\mathbf{o}\in\mathcal{O}} \frac{\beta^\mathbf{o}(\bm{\mu})}{P_{1}^{\mathbf{o}} } \leq \sup_{\bm{\mu}\in D} \sum_{\mathbf{o}\in\mathcal{O}} \frac{\alpha^\mathbf{o}(\bm{\mu})}{P_{0}^{\mathbf{o}}}
\end{align}
by showing 
\begin{align} \label{eq:stepopt_pf_6}
& \lim_{m\to \infty} \sup_{\bm{\mu}\in D} \sum_{\mathbf{o}\in\mathcal{O}} \frac{\alpha_{m}^\mathbf{o}(\bm{\mu})}{P_{0,m}^{\mathbf{o}}} = \sup_{\bm{\mu}\in D} \sum_{\mathbf{o}\in\mathcal{O}} \frac{\alpha^\mathbf{o}(\bm{\mu})}{P_{0}^{\mathbf{o}}} \\ \label{eq:stepopt_pf_7}
& \lim_{m\to \infty} \inf_{\bm{\mu}\in D} \sum_{\mathbf{o}\in\mathcal{O}} \frac{\beta_{m}^\mathbf{o}(\bm{\mu})}{P_{1,m}^{\mathbf{o}} } = \inf_{\bm{\mu}\in D} \sum_{\mathbf{o}\in\mathcal{O}} \frac{\beta^\mathbf{o}(\bm{\mu})}{P_{1}^{\mathbf{o}}}.
\end{align}

\textbf{Proving \eqref{eq:stepopt_pf_6}.}
Because $P_{0,m}^{\mathbf{o}} \uparrow P_{0}^{\mathbf{o}}$ not depending on $\bm{\mu}$, $P_{0}^{\mathbf{o}}$ bounded away from $0$, and $\alpha_m^\mathbf{o}(\bm{\mu})\leq 1$, we have that for any $\epsilon_0>0$, there exists $m_0$ such that $\forall~m>m_0, \bm{\mu}\in D$, 
\begin{align*}
\sum_{\mathbf{o}\in\mathcal{O}} \frac{\alpha_m^\mathbf{o}(\bm{\mu})}{P_{0}^{\mathbf{o}}}  \leq \sum_{\mathbf{o}\in\mathcal{O}} \frac{\alpha_{m}^\mathbf{o}(\bm{\mu})}{P_{0,m}^{\mathbf{o}}} < \sum_{\mathbf{o}\in\mathcal{O}} \frac{\alpha_m^\mathbf{o}(\bm{\mu})}{P_{0}^{\mathbf{o}}} + \epsilon_0,
\end{align*}
giving 
\begin{align} \label{eq:stepopt_pf_8}
\lim_{m\to \infty} \sup_{\bm{\mu}\in D} \sum_{\mathbf{o}\in\mathcal{O}} \frac{\alpha_{m}^\mathbf{o}(\bm{\mu})}{P_{0,m}^{\mathbf{o}}} = \lim_{m\to \infty} \sup_{\bm{\mu}\in D} \sum_{\mathbf{o}\in\mathcal{O}} \frac{\alpha_{m}^\mathbf{o}(\bm{\mu})}{P_{0}^{\mathbf{o}}}.
\end{align}
Next, as $\alpha_{m}^\mathbf{o}(\bm{\mu})\uparrow\alpha^\mathbf{o}(\bm{\mu})$, we have $\sum_{\mathbf{o}\in\mathcal{O}} \frac{\alpha_{m}^\mathbf{o}(\bm{\mu})}{P_{0}^{\mathbf{o}}} \uparrow \sum_{\mathbf{o}\in\mathcal{O}} \frac{\alpha^\mathbf{o}(\bm{\mu})}{P_{0}^{\mathbf{o}}}$, giving
\begin{align*}
\lim_{m \to \infty} \sup_{\bm{\mu}\in D} \sum_{\mathbf{o}\in\mathcal{O}} \frac{\alpha_{m}^\mathbf{o}(\bm{\mu})}{P_{0}^{\mathbf{o}}} \leq \sup_{\bm{\mu}\in D} \sum_{\mathbf{o}\in\mathcal{O}} \frac{\alpha^\mathbf{o}(\bm{\mu})}{P_{0}^{\mathbf{o}}}.
\end{align*}
On the other hand, for any $\epsilon_1>0$, let 
\begin{align*}
& E=\{\bm{\mu}:\sum_{\mathbf{o}\in\mathcal{O}} \frac{\alpha^\mathbf{o}(\bm{\mu})}{P_{0}^{\mathbf{o}}}>\sup_{\bm{\mu}\in D}\sum_{\mathbf{o}\in\mathcal{O}} \frac{\alpha^\mathbf{o}(\bm{\mu})}{P_{0}^{\mathbf{o}}}-\epsilon_1\} \\
& E_{m}=\{\bm{\mu}:\sum_{\mathbf{o}\in\mathcal{O}} \frac{\alpha_{m}^\mathbf{o}(\bm{\mu})}{P_{0}^{\mathbf{o}}}>\sup_{\bm{\mu}\in D}\sum_{\mathbf{o}\in\mathcal{O}} \frac{\alpha^\mathbf{o}(\bm{\mu})}{P_{0}^{\mathbf{o}}}-\epsilon_1\}.
\end{align*}
Notice that $\sum_{\mathbf{o}\in\mathcal{O}} \frac{\alpha_{m}^\mathbf{o}(\bm{\mu})}{P_{0}^{\mathbf{o}}} \uparrow \sum_{\mathbf{o}\in\mathcal{O}} \frac{\alpha^\mathbf{o}(\bm{\mu})}{P_{0}^{\mathbf{o}}}$, we have $E_{m}\subset E_{m+1}$, $\forall~n$, and $E=\bigcup_{m}E_{m}$.
Therefore, for any $\epsilon_1>0$, $\exists ~m_{1}$ such that for all $m>m_{1},$ $E_{m}\neq\emptyset$, which implies 
\begin{align*}
\lim_{m \to \infty} \sup_{\bm{\mu}\in D} \sum_{\mathbf{o}\in\mathcal{O}} \frac{\alpha_{m}^\mathbf{o}(\bm{\mu})}{P_{0}^{\mathbf{o}}} \geq \sup_{\bm{\mu}\in D} \sum_{\mathbf{o}\in\mathcal{O}} \frac{\alpha^\mathbf{o}(\bm{\mu})}{P_{0}^{\mathbf{o}}}.
\end{align*}
Therefore,
\begin{align}\label{eq:stepopt_pf_9}
\lim_{m\to \infty} \sup_{\bm{\mu}\in D} \sum_{\mathbf{o}\in\mathcal{O}} \frac{\alpha_{m}^\mathbf{o}(\bm{\mu})}{P_{0}^{\mathbf{o}}} = \sup_{\bm{\mu}\in D} \sum_{\mathbf{o}\in\mathcal{O}} \frac{\alpha^\mathbf{o}(\bm{\mu})}{P_{0}^{\mathbf{o}}}.
\end{align}
Combining \eqref{eq:stepopt_pf_8} and \eqref{eq:stepopt_pf_9} to have \eqref{eq:stepopt_pf_6}. 

\textbf{Proving \eqref{eq:stepopt_pf_7}.}
Due to the same reason of \eqref{eq:stepopt_pf_8} we have 
\begin{align}
\label{eq:stepopt_pf_10}
\lim_{m\to \infty} \inf_{\bm{\mu}\in D} \sum_{\mathbf{o}\in\mathcal{O}} \frac{\beta_{m}^\mathbf{o}(\bm{\mu})}{P_{1,m}^{\mathbf{o}}} = \lim_{m\to \infty} \inf_{\bm{\mu}\in D} \sum_{\mathbf{o}\in\mathcal{O}} \frac{\beta_{m}^\mathbf{o}(\bm{\mu})}{P_{1}^{\mathbf{o}}}.
\end{align}
Since $\beta_{m}(\bm{\mu}) \uparrow \beta(\bm{\mu})$, we have $\sum_{\mathbf{o}\in\mathcal{O}} \frac{\beta_{m}^\mathbf{o}(\bm{\mu})}{P_{1}^{\mathbf{o}} } \uparrow \sum_{\mathbf{o}\in\mathcal{O}} \frac{\beta^\mathbf{o}(\bm{\mu})}{P_{1}^{\mathbf{o}}}$, giving 
\begin{align}\label{eq:stepopt_pf_13}
\lim_{m\to\infty} \inf_{\bm{\mu}\in D} \sum_{\mathbf{o}\in\mathcal{O}} \frac{\beta_{m}^\mathbf{o}(\bm{\mu})}{P_{1}^{\mathbf{o}} } \leq \inf_{\bm{\mu}\in D} \sum_{\mathbf{o}\in\mathcal{O}} \frac{\beta^\mathbf{o}(\bm{\mu})}{P_{1}^{\mathbf{o}}}.
\end{align}
Define $c_{m}=\sum_{\mathbf{o}\in\mathcal{O}}\sum_{\mathbf{b}<-m~or~\mathbf{b}>m}\frac{P_{1}(I_{\mathbf{b}}^{\mathbf{o}})}{P_{1}^\mathbf{o}}$.
Then $c_{m}\downarrow 0$ and $\forall \bm{\mu}\in D$, $\sum_{\mathbf{o}\in\mathcal{O}} \frac{\beta_{m}^\mathbf{o}(\bm{\mu})}{P_{1}^{\mathbf{o}} }+c_{m}\downarrow \sum_{\mathbf{o}\in\mathcal{O}} \frac{\beta^\mathbf{o}(\bm{\mu})}{P_{1}^{\mathbf{o}}}$.
Furthermore, for any $\epsilon_1>0$, let 
\begin{align*}
& E=\{\bm{\mu}:\sum_{\mathbf{o}\in\mathcal{O}} \frac{\beta^\mathbf{o}(\bm{\mu})}{P_{1}^{\mathbf{o}}}<\lim_{m\to \infty}\inf_{\mu\in D}\left(\sum_{\mathbf{o}\in\mathcal{O}} \frac{\beta_{m}^\mathbf{o}(\bm{\mu})}{P_{1}^{\mathbf{o}} }+c_{m}\right)+\epsilon_1\}\\
& E_{m}=\{\mu:\sum_{\mathbf{o}\in\mathcal{O}} \frac{\beta_{m}^\mathbf{o}(\bm{\mu})}{P_{1}^{\mathbf{o}} }+c_{m}<\lim_{m\to \infty}\inf_{\mu\in D}\left(\sum_{\mathbf{o}\in\mathcal{O}} \frac{\beta_{m}^\mathbf{o}(\bm{\mu})}{P_{1}^{\mathbf{o}} }+c_{m}\right)+\epsilon_1\}.
\end{align*}
Since $\sum_{\mathbf{o}\in\mathcal{O}} \frac{\beta_{m}^\mathbf{o}(\bm{\mu})}{P_{1}^{\mathbf{o}} }+c_{m}\downarrow \sum_{\mathbf{o}\in\mathcal{O}} \frac{\beta^\mathbf{o}(\bm{\mu})}{P_{1}^{\mathbf{o}} }$, we have $E_{m}\subset E_{m+1}$,
$\forall~m$, and $E=\bigcup_{m}E_{m}$.  
For any $\epsilon_1>0$, since for any $m$, $E_m \neq \emptyset$, we have $E\neq\emptyset$ and hence 
\begin{align}\label{eq:stepopt_pf_14}
\inf_{\bm{\mu}\in D} \sum_{\mathbf{o}\in\mathcal{O}} \frac{\beta^\mathbf{o}(\bm{\mu})}{P_{1}^{\mathbf{o}}} \leq\lim_{m\to\infty} \inf_{\bm{\mu}\in D} \left(\sum_{\mathbf{o}\in\mathcal{O}} \frac{\beta_{m}^\mathbf{o}(\bm{\mu})}{P_{1}^{\mathbf{o}} }+c_{m} \right)=\lim_{m\to\infty} \inf_{\bm{\mu}\in D} \sum_{\mathbf{o}\in\mathcal{O}} \frac{\beta_{m}^\mathbf{o}(\bm{\mu})}{P_{1}^{\mathbf{o}}}.
\end{align}
Combining \eqref{eq:stepopt_pf_13} and \eqref{eq:stepopt_pf_14} we obtain 
\begin{align}\label{eq:stepopt_pf_11}
\inf_{\bm{\mu}\in D} \sum_{\mathbf{o}\in\mathcal{O}} \frac{\beta^\mathbf{o}(\bm{\mu})}{P_{1}^{\mathbf{o}}} =\lim_{m\to\infty} \inf_{\bm{\mu}\in D} \sum_{\mathbf{o}\in\mathcal{O}} \frac{\beta_{m}^\mathbf{o}(\bm{\mu})}{P_{1}^{\mathbf{o}}}.
\end{align}
Combining \eqref{eq:stepopt_pf_10} and \eqref{eq:stepopt_pf_11} to have \eqref{eq:stepopt_pf_7}. 
Then finally we proved \eqref{eq:stepopt_pf_5}. 

\textbf{Finally we show the power of the sign test upper bounds that of any test in $\mathcal{S}(\omega)$ using \eqref{eq:stepopt_pf_5}.}
\eqref{eq:stepopt_pf_5} further gives 
\begin{align} \label{eq:stepopt_pf_12}
\inf_{\bm{\mu}>0} \sum_{\mathbf{o}\in\mathcal{O}} \frac{\beta^\mathbf{o}(\bm{\mu})}{P_{1}^{\mathbf{o}} }  \leq \inf_{\bm{\mu}\in D} \sum_{\mathbf{o}\in\mathcal{O}} \frac{\beta^\mathbf{o}(\bm{\mu})}{P_{1}^{\mathbf{o}} } \leq \sup_{\bm{\mu}\in D} \sum_{\mathbf{o}\in\mathcal{O}} \frac{\alpha^\mathbf{o}(\bm{\mu})}{P_{0}^{\mathbf{o}}} \leq \sup_{\bm{\mu} > 0} \sum_{\mathbf{o}\in\mathcal{O}} \frac{\alpha^\mathbf{o}(\bm{\mu})}{P_{0}^{\mathbf{o}}} 
\end{align}
Recall that for every orthant $\mathbf{o}$, $P^{\mathbf{o}}_0 = 2^{-n}$.
Multiplying both sides of \eqref{eq:stepopt_pf_12} by $2^{-n}$ we have 
\begin{align}
\inf_{\bm{\mu}>0} 2^{-n} \sum_{\mathbf{o}\in\mathcal{O}} \frac{\beta^\mathbf{o}(\bm{\mu})}{P_{1}^{\mathbf{o}} } \leq \sup_{\bm{\mu} > 0} \sum_{\mathbf{o}\in\mathcal{O}} \alpha^\mathbf{o}(\bm{\mu}) = \sup_{\bm{\mu} > 0} \alpha(\bm{\mu}) \leq \alpha,
\end{align}
where $\alpha$ is the size of the test. 
For any $\epsilon>0$, there exists an $\bm{\mu}'$ such that $2^{-n} \sum_{\mathbf{o}\in\mathcal{O}} \frac{\beta^\mathbf{o}(\bm{\mu}')}{P_{1}^{\mathbf{o}} } \leq \alpha+\epsilon$.
For this specific $\bm{\mu}'$, the power for each orthant is weighted by $1/P_1^o$ in the upper bound. 
To maximize the overall power $\beta(\bm{\mu}')=\sum_{\mathbf{o}\in\mathcal{O}} \beta^\mathbf{o}(\bm{\mu}')$, we start from the orthant with the largest $P_1^\mathbf{o}$, maximizing its power $\beta^\mathbf{o}(\bm{\mu}')$ by letting the test $\phi$ to be $1$ for that orthant. 
Then we go to the second largest and keep doing it until the inequality becomes equal. 
This is indeed the sign test, giving that $\beta(\bm{\mu}')$ is no larger than the power of the sign test. 
Together with the fact that $\beta(\bm{\mu}')$ serves as an upper bound for the worse-case power of $\phi$, it implies that $\phi^{S}$ is maximin among all tests in ${\cal S}(\omega)$. 
Furthermore, for any $\phi\in{\cal S}$, $\exists\;\omega$ such that $\phi\in{\cal S}(\omega)$. 
By noting that $\phi^S$ is a maximin $\alpha$-level test in any ${\cal S}(\omega)$ we complete the proof. 
\end{proof}

\subsection{Proof of Lemma \ref{lem:sim_apr}}
\begin{proof}
Consider any $\phi(\mathbf{Y}) \in \mathcal{B}_n$.
For any $\omega>0$,
define $\psi_{(\omega)}=\sum_{\mathbf{o}\in\mathcal{O}}\sum_{-\infty<\mathbf{b}<\infty} \psi^{\mathbf{o}}_\mathbf{b} \ind_{I^\mathbf{o}_\mathbf{b}}$, where $\psi_{\mathbf{b}}^{\mathbf{o}}=\frac{1}{\vert I^\mathbf{o}_\mathbf{b}\vert}\int_{I^\mathbf{o}_\mathbf{b}}\phi(\mathbf{y})$ and $\vert I^\mathbf{o}_\mathbf{b}\vert$ is the volume of the box $I^\mathbf{o}_\mathbf{b}$.
We next show that there exists a small enough $\omega$ such that $\phi\doteqdel{\epsilon}\psi_{(\omega)}$. 

Notice that by fixing $\omega\in\mathbb{R}$, any real positive vector $\bm{\mu}=(\mu_1,\cdots,\mu_n)$ can be written as $\mu_j=(1+\omega)^{d_j}(1+\omega'_j)$, for $\mathbf{d}=(d_1,\cdots,d_n)$ and $\bm{\omega}'=(\omega'_1,\cdots,\omega'_n)$, where $0\leq \bm{\omega}' <\omega$ and we recall that the inequality of the vector $\bm{\omega}'$ is element-wise. 
Similar to Step 2 in the proof sketch of Theorem \ref{thm: opt_signtest}, define $\tilde{I}_{b}^{+}(\omega'_j)=\frac{1}{1+\omega'_j}I_{b}^{+}=(\frac{(1+\omega)^{b}}{1+\omega'},\frac{(1+\omega)^{b+1}}{1+\omega'_j}]$
and $\tilde{I}_{b}^{-}(\omega'_j)=\frac{1}{1+\omega'_j}I_{b}^{-}(\omega'_j)=[-\frac{(1+\omega)^{b+1}}{1+\omega'_j},-\frac{(1+\omega)^{b}}{1+\omega'_j}$).
Then we can define the rescaled box $\tilde{I}^\mathbf{o}_\mathbf{b} = \tilde{I}^{o_1}_{b_1}(\omega'_1) \times \cdots \times \tilde{I}^{o_n}_{b_n}(\omega'_n)$. 
Note that the element-wise multiplication by $\bm{\mu}$ maps the rescaled box $\tilde{I}^\mathbf{o}_\mathbf{b}$ to the original box $I^\mathbf{o}_{\mathbf{b}+\mathbf{d}}$.

First, for any index $\mathbf{b}$, we have
\begin{align*}
 & \int_{\tilde{I}^\mathbf{o}_\mathbf{b}}\left[\phi(\bm{\mu} \circ \mathbf{z})-\psi(\bm{\mu} \circ \mathbf{z})\right]d\mathbf{z}=\int_{\tilde{I}^\mathbf{o}_\mathbf{b}}\phi(\bm{\mu} \circ \mathbf{z})d\mathbf{z}-\vert\tilde{I}_{\mathbf{b}}^{\mathbf{o}}\vert\psi_{\mathbf{b}+\mathbf{d}}^{\mathbf{o}}\\
& =  \frac{1}{\prod_j \mu_j}\int_{I_{\mathbf{b}+\mathbf{d}}^{\mathbf{o}}}\phi(\mathbf{z})d\mathbf{z}-\vert\tilde{I}_{\mathbf{b}}^{\mathbf{o}}\vert\psi_{\mathbf{b}+\mathbf{d}}^{\mathbf{o}}
= \frac{1}{\prod_j\mu_j}\vert I_{\mathbf{b}+\mathbf{d}}^{\mathbf{o}} \vert\psi_{\mathbf{b}+\mathbf{d}}^{\mathbf{o}}-\vert\tilde{I}_{\mathbf{b}}^{\mathbf{o}}\vert\psi_{\mathbf{b}+\mathbf{d}}^{\mathbf{o}}=0.
\end{align*}

Let us define $\underline{f}_{\mathbf{b}}^{\mathbf{o}}(\bm{\omega}')=\inf_{\mathbf{z}\in\tilde{I}_{\mathbf{b}}^{\mathbf{o}}}f(\mathbf{z})$.
For any $\bm{\mu}$ specified by $\mathbf{d}$ and $\bm{\omega}'$, the difference of the sizes of the two tests $\phi$ and $\psi$ can be upper bounded as
\begin{align*}
& \vert\E_{0}[\phi(\mathbf{Y})-\psi(\mathbf{Y})]\vert = \bigg\vert\int_{\mathbb{R}^n}\left[\phi(\bm{\mu}\circ \mathbf{z})-\psi(\bm{\mu}\circ \mathbf{z})\right]f(\mathbf{z})d\mathbf{z}\bigg\vert \\
& = \bigg\vert \sum_{\mathbf{o}\in\mathcal{O}}\sum_{-\infty<\mathbf{b}<\infty}\int_{\tilde{I}_{\mathbf{b}}^{\mathbf{o}}}\left[\phi(\bm{\mu}\circ \mathbf{z})-\psi(\bm{\mu}\circ \mathbf{z})\right]f(\mathbf{z})d\mathbf{z} \bigg\vert\\
& = \bigg\vert \sum_{\mathbf{o}\in\mathcal{O}}\sum_{-\infty<\mathbf{b}<\infty}\int_{\tilde{I}_{\mathbf{b}}^{\mathbf{o}}}\left[\phi(\bm{\mu}\circ \mathbf{z})-\psi(\bm{\mu}\circ \mathbf{z})\right]\left[f(\mathbf{z})-\underline{f}_{\mathbf{b}}^{\mathbf{o}}(\bm{\omega}')\right]d\mathbf{z} \bigg\vert \\
& \leq \sum_{\mathbf{o}\in\mathcal{O}} \sum_{-\infty < \mathbf{b}<\infty} \int_{\tilde{I}_{\mathbf{b}}^{\mathbf{o}}} \left[f(\mathbf{z})-\underline{f}_{\mathbf{b}}^{\mathbf{o}}(\bm{\omega}')\right]d\mathbf{z}\\
& \leq \sum_{\mathbf{o}\in\mathcal{O}} \left[ \sum_{-m \leq \mathbf{b}\leq m} \int_{\tilde{I}_{\mathbf{b}}^{\mathbf{o}}} \left[f(\mathbf{z})-\underline{f}_{\mathbf{b}}^{\mathbf{o}}(\bm{\omega}')\right]d\mathbf{z} +\sum_{\mathbf{b}: \exists \vert b_j\vert >m} \int_{\tilde{I}_{\mathbf{b}}^{\mathbf{o}}} \left[f(\mathbf{z})-\underline{f}_{\mathbf{b}}^{\mathbf{o}}(\bm{\omega}')\right]d\mathbf{z} \right]
\end{align*}
As $\omega\to 0$ and $m\to\infty$, the above expression goes uniformly
to zero for all $\forall~\bm{\omega}'<\omega$. Therefore $\exists~\omega_{0},m_{0}$,
such that $\forall\;\omega<\omega_{0}$, $m>m_{0}$, we have 
\begin{align*}
& \sum_{\mathbf{o}\in\mathcal{O}} \sum_{-m \leq \mathbf{b}\leq m} \int_{\tilde{I}_{\mathbf{b}}^{\mathbf{o}}} \left[f(\mathbf{z})-\underline{f}_{\mathbf{b}}^{\mathbf{o}}(\bm{\omega}')\right]d\mathbf{z} \leq \frac{\epsilon}{2}, ~~~~\forall \bm{\omega}'<\omega,\\
& \sum_{\mathbf{o}\in\mathcal{O}} \sum_{\mathbf{b}: \exists \vert b_j\vert >m} \int_{\tilde{I}_{\mathbf{b}}^{\mathbf{o}}} \left[f(\mathbf{z})-\underline{f}_{\mathbf{b}}^{\mathbf{o}}(\bm{\omega}')\right]d\mathbf{z} \leq \frac{\epsilon}{2}.
\end{align*}
Thus $\forall \mu>0$, $\vert\E_{0}[\phi(\mathbf{Y})-\psi(\mathbf{Y})]\vert\leq\epsilon$, which implies that the sizes of the two tests are within $\epsilon$-distance of each other. 

Similarly, we can bound the difference between the power of the two tests. Let $\underline{\underline{f_{\mathbf{b}}^{\mathbf{o}}}}(\bm{\omega}')=\inf_{\mathbf{z}\in\tilde{I}_{\mathbf{b}}^{\mathbf{o}}}f(\mathbf{z}-\delta)$.
For any $\bm{\mu}$ specified by $\mathbf{d}$ and $\bm{\omega}'$, the difference of the power
of $\phi$ and $\psi$ can be written as
\begin{align*}
& \vert\E_{1}[\phi(\mathbf{Y})-\psi(\mathbf{Y})]\vert = \bigg\vert\int_{\mathbb{R}^n}\left[\phi(\bm{\mu}\circ \mathbf{z})-\psi(\bm{\mu}\circ \mathbf{z})\right]f(\mathbf{z}-\delta)d\mathbf{z}\bigg\vert \\
& = \bigg\vert \sum_{\mathbf{o}\in\mathcal{O}}\sum_{-\infty<\mathbf{b}<\infty}\int_{\tilde{I}_{\mathbf{b}}^{\mathbf{o}}}\left[\phi(\bm{\mu}\circ \mathbf{z})-\psi(\bm{\mu}\circ \mathbf{z})\right]f(\mathbf{z}-\delta)d\mathbf{z} \bigg\vert\\
& = \bigg\vert \sum_{\mathbf{o}\in\mathcal{O}}\sum_{-\infty<\mathbf{b}<\infty}\int_{\tilde{I}_{\mathbf{b}}^{\mathbf{o}}}\left[\phi(\bm{\mu}\circ \mathbf{z})-\psi(\bm{\mu}\circ \mathbf{z})\right]\left[f(\mathbf{z}-\delta)-\underline{\underline{f_{\mathbf{b}}^{\mathbf{o}}}}(\bm{\omega}')\right]d\mathbf{z} \bigg\vert \\
& \leq \sum_{\mathbf{o}\in\mathcal{O}} \sum_{-\infty < \mathbf{b}<\infty} \int_{\tilde{I}_{\mathbf{b}}^{\mathbf{o}}} \left[f(\mathbf{z}-\delta)-\underline{\underline{f_{\mathbf{b}}^{\mathbf{o}}}}(\bm{\omega}')\right]d\mathbf{z}\\
& \leq \sum_{\mathbf{o}\in\mathcal{O}} \left[ \sum_{-m \leq \mathbf{b}\leq m} \int_{\tilde{I}_{\mathbf{b}}^{\mathbf{o}}} \left[f(\mathbf{z}-\delta)-\underline{\underline{f_{\mathbf{b}}^{\mathbf{o}}}}(\bm{\omega}')\right]d\mathbf{z} +\sum_{\mathbf{b}: \exists \vert b_j\vert >m} \int_{\tilde{I}_{\mathbf{b}}^{\mathbf{o}}} \left[f(\mathbf{z}-\delta)-\underline{\underline{f_{\mathbf{b}}^{\mathbf{o}}}}(\bm{\omega}')\right]d\mathbf{z} \right]\end{align*}

Similar to the argument in the previous case, as $\omega\to 0$ and $m\to\infty$, the above expression
goes uniformly to zero for all $\bm{\omega}'<\omega$. 
Hence there exists $\omega_{1},m_{1}$, such that $\forall\;\omega<\omega_{1}$, $m>m_{1}$, 
\begin{align*}
&\sum_{\mathbf{o}\in\mathcal{O}} \sum_{-m \leq \mathbf{b}\leq m} \int_{\tilde{I}_{\mathbf{b}}^{\mathbf{o}}} \left[f(\mathbf{z}-\delta)-\underline{\underline{f_{\mathbf{b}}^{\mathbf{o}}}}(\bm{\omega}')\right]d\mathbf{z}<\frac{\epsilon}{2}, ~~~~\forall \bm{\omega}'<\omega, \\
& \sum_{\mathbf{o}\in\mathcal{O}}\sum_{\mathbf{b}: \exists \vert b_j\vert >m} \int_{\tilde{I}_{\mathbf{b}}^{\mathbf{o}}} \left[f(\mathbf{z}-\delta)-\underline{\underline{f_{\mathbf{b}}^{\mathbf{o}}}}(\bm{\omega}')\right] d\mathbf{z} < \frac{\epsilon}{2}.
\end{align*}

Thus $\forall \bm{\mu}>0$, $\vert\E_{1}[\phi(\mathbf{Y})-\psi(\mathbf{Y})]\vert\leq\epsilon$ which implies the distance between the power of the two tests is bounded by $\epsilon$. 
Therefore, for any $\epsilon>0,$ by selecting $\omega<\min(\omega_{0},\omega_{1})$ and $m>\max(m_{0},m_{1})$, we have a test $\psi_{(\omega)}\in{\cal S}$, such that $\forall\;\mu>0$, $\vert\E_{0}[\phi(\mathbf{Y})-\psi(\mathbf{Y})]\vert\leq\epsilon$ and $\vert\E_{1}[\phi(\mathbf{Y})-\psi(\mathbf{Y})]\vert\leq\epsilon$. 
As a result, $\phi\doteqdel{\epsilon}\psi_{(\omega)}$.
\end{proof}

\end{document}